\documentclass{article}
\usepackage{color}
\usepackage{microtype}
\usepackage{graphicx}
\usepackage{subfigure}
\usepackage{booktabs} % for professional tables
\usepackage{amsmath}
\usepackage{amssymb}
\usepackage{bbm}
\usepackage[T1]{fontenc}
\usepackage{arydshln}
\makeatletter
\renewcommand{\maketag@@@}[1]{\hbox{\m@th\normalsize\normalfont#1}}%
\makeatother
\makeatletter
\def\thickhline{%
	\noalign{\ifnum0=`}\fi\hrule \@height \thickarrayrulewidth \futurelet
	\reserved@a\@xthickhline}
\def\@xthickhline{\ifx\reserved@a\thickhline
	\vskip\doublerulesep
	\vskip-\thickarrayrulewidth
	\fi
	\ifnum0=`{\fi}}
\makeatother

\newlength{\thickarrayrulewidth}
\usepackage{multirow}
\usepackage{array}
\newcolumntype{P}[1]{>{\centering\arraybackslash}p{#1}}
\newcolumntype{M}[1]{>{\centering\arraybackslash}m{#1}}

\usepackage{amsmath,amsfonts,bm,amssymb,mathtools}

\def\vzero{{\bm{0}}}

\def\vtheta{{\bm{\theta}}}

\def\vf{{\bm{f}}}
\def\vg{{\bm{g}}}

\def\vq{{\bm{q}}}

\def\vs{{\bm{s}}}

\def\vu{{\bm{u}}}

\def\vx{{\bm{x}}}
\def\vy{{\bm{y}}}
\def\vz{{\bm{z}}}

\def\vtheta{{\bm{\theta}}}

\DeclareMathAlphabet{\mathsfit}{\encodingdefault}{\sfdefault}{m}{sl}
\SetMathAlphabet{\mathsfit}{bold}{\encodingdefault}{\sfdefault}{bx}{n}

\usepackage{caption}
\usepackage[noend]{algorithmic}

\usepackage{makecell}
\usepackage{ntheorem}

\newtheorem{thm}{Theorem}[section]

\newtheorem{ex}[thm]{Example}

\newtheorem{remark}[thm]{Remark}
\newtheorem{definition}{Definition}[section]
\usepackage{hyperref}

\usepackage[most]{tcolorbox}
\definecolor{myfavblue}{rgb}{0.05, 0.2, 0.8}
\definecolor{black}{rgb}{0.0, 0.0, 0.0}
\definecolor{wine}{rgb}{0.5333333333333333, 0.13333333333333333, 0.3333333333333333}

\usepackage[accepted]{icml2024}

\usepackage[textsize=tiny]{todonotes}
\usepackage[textsize=tiny]{todonotes}

\usepackage[capitalize,noabbrev]{cleveref}

\icmltitlerunning{Submission and Formatting Instructions for ICML 2024}
\icmltitlerunning{FESSNC: Fast Exponentially Stable and Safe Neural Controller}
\begin{document}

\twocolumn[
\icmltitle{FESSNC: Fast Exponentially Stable and Safe Neural Controller}

% It is OKAY to include author information, even for blind
% submissions: the style file will automatically remove it for you
% unless you've provided the [accepted] option to the icml2021
% package.

% List of affiliations: The first argument should be a (short)
% identifier you will use later to specify author affiliations
% Academic affiliations should list Department, University, City, Region, Country
% Industry affiliations should list Company, City, Region, Country

% You can specify symbols, otherwise they are numbered in order.
% Ideally, you should not use this facility. Affiliations will be numbered
% in order of appearance and this is the preferred way.
\icmlsetsymbol{equal}{*}

\begin{icmlauthorlist}
\icmlauthor{Jingdong Zhang}{equal,MS,IICS}
\icmlauthor{Luan Yang}{equal,IICS}
\icmlauthor{Qunxi Zhu}{IICS,MOE,AIlab}
\icmlauthor{Wei Lin}{MS,IICS,MOE,AIlab}
\end{icmlauthorlist}

\icmlaffiliation{MS}{School of Mathematical Sciences, LMNS, and SCMS, Fudan University, China.}
\icmlaffiliation{IICS}{Research Institute of Intelligent Complex Systems, Fudan University, China.}
\icmlaffiliation{MOE}{State Key Laboratory of Medical Neurobiology and MOE Frontiers Center for Brain Science, Institutes of Brain Science, Fudan University, China.}
\icmlaffiliation{AIlab}{Shanghai Artificial Intelligence Laboratory, China}

\icmlcorrespondingauthor{Qunxi Zhu}{qxzhu16@fudan.edu.cn}
\icmlcorrespondingauthor{Wei Lin}{wlin@fudan.edu.cn}
% You may provide any keywords that you
% find helpful for describing your paper; these are used to populate
% the "keywords" metadata in the PDF but will not be shown in the document
\icmlkeywords{Machine Learning, ICML}

\vskip 0.3in
]

% this must go after the closing bracket ] following \twocolumn[ ...

% This command actually creates the footnote in the first column
% listing the affiliations and the copyright notice.
% The command takes one argument, which is text to display at the start of the footnote.
% The \icmlEqualContribution command is standard text for equal contribution.
% Remove it (just {}) if you do not need this facility.

%\printAffiliationsAndNotice{}  % leave blank if no need to mention equal contribution
\printAffiliationsAndNotice{\icmlEqualContribution} % otherwise use the standard text.

\begin{abstract}
   In order to stabilize nonlinear systems modeled by stochastic differential equations, we design a \textbf{F}ast \textbf{E}xponentially \textbf{S}table and \textbf{S}afe \textbf{N}eural \textbf{C}ontroller (FESSNC) for \textit{fast} learning controllers. Our framework is parameterized by neural networks, and realizing both rigorous exponential stability and safety guarantees.  Concretely, we design heuristic methods to learn the exponentially stable and the safe controllers, respectively, in light of the classic stochastic exponential stability theory and our established theorem on guaranteeing the almost-sure safety for stochastic dynamics.  More significantly, to rigorously ensure the stability and the safety guarantees for the learned controllers, we develop a projection operator, projecting to the space of exponentially-stable and safe controllers. To reduce the high computation cost of solving the projection operation, approximate projection operators are delicately proposed with closed forms that map the learned controllers to the target controller space. Furthermore, we employ Hutchinson's trace estimator for a scalable unbiased estimate of the Hessian matrix that is used in the projection operator, which thus allows for computation cost reduction and therefore can accelerate the training and testing processes. {\color{black}More importantly, our approximate projection operations can be applied to the nonparametric control methods to improve their stability and safety performance.} We empirically demonstrate the superiority of the FESSNC over the existing methods.
\end{abstract}

\section{Introduction}
    Stabilizing the nonlinear systems modeled by stochastic differential equations (SDEs) is a challenging focal task in many fields of mathematics and engineering. The neural networks (NNs) equipped with the suitable loss of stability conditions have remarkable abilities in learning the stabilization controllers for the underlying systems, such as the neural Lyapunov control for ordinary differential equations (ODEs)~\citep{chang2020neural} and the neural stochastic control for SDEs~\citep{zhang2022neural}. However, these frameworks, which heuristically train the NNs on the finite samples from the data space, are suffered from the lack of a rigorous stability guarantee for the learned neural controllers, i.e., satisfying the anticipated conditions on the whole data space.

    Safety is another major concern for controlled dynamics, such as robotic and automotive systems~\citep{wang2016multi}. Recent advances in the control barrier function theory pave a direct way to impose the safety for the nonlinear SDEs.  Several algorithms have theoretically investigated the structures of the barrier functions and designed the optimal control via quadratic programming (QP) or Sum-of-Squares (SOS) optimization according to the theoretical results~\citep{sarkar2020high,mazouzsafety}. Although all these algorithms try to obtain the controllers that rigorously satisfy the safety conditions, they either suffer high computational costs or only consider the safety with probability $\delta<1$ rather than the  probability one which is urgently required in the safety-critical applications.

    In this paper, our goal is to provide a framework for fast learning the neural controllers with rigorous exponential stability and safety guarantees such that the trajectories of the controlled SDEs exponentially converge to the target equilibria and remain in the safe region all the time. %To this end, we propose the fast exponentially stable and safe neural controller (FESSNC), a learning-based framework that combines the projection operator to the space of exponentially stable and safe controllers, and the accelerating tricks for trace estimation. 
    The major contributions of this paper are summarized as follows.
\begin{itemize}
    \item We train the neural controllers with the exponential stability conditions, and impose an exponential stability guarantee for the learned controller by projecting it into the space of the exponentially stable controllers. To avoid the high computational complexity of solving the projection operation, we propose an approximate projection operator $\hat{\pi}_{es}$ with a closed form that can map the controller to the target space if and only if the state space of the controlled SDEs is bounded.
    \item %To solve the problem that exponentially stable SDEs can have unbounded state space and cannot remain in the specified safe region,  
    We provide a theorem for restricting the controlled SDEs in the predefined bounded safe region based on the zeroing barrier functions. Similarly, we train the neural controller according to the safety loss and project it to the space of safe controllers with the approximate projection operator $\hat{\pi}_{sf}$. Theoretically, it can be verified that one can combine $\hat{\pi}_{es}$ and $\hat{\pi}_{sf}$ to provide both exponential stability and safety guarantees for the neural controllers.
    \item We introduce an unbiased stochastic estimator of trace estimation of the Hessian matrix involved in the projection operator with cheap computation.
    \item We demonstrate the efficacy of the FESSNC on a variety of classic control problems, such as, the synchronization of coupled oscillators. {\color{black}In addition to the neural network controllers, we also show our framework can be applied to the nonparametric controllers to provide stability and safety guarantee for them.} Our code is available at  \href{https://github.com/jingddong-zhang/FESSNC}{\texttt{github.com/jingddong-zhang/FESSNC}}.
\end{itemize}

\section{Preliminaries}
\subsection{Problem Formulation}

To begin with, we consider the SDE in a general form of
\begin{equation}\label{SDE0}
	\mathrm{d}{\vx}(t) = f({\vx}(t))\mathrm{d}t + g({\vx}(t))\mathrm{d}B_t,
	~t\ge0,
\end{equation}
where $\vx(0)=\vx_0\in\mathbb{R}^d,~\{\vx(t)\}_{t\ge0}\subset\mathcal{X}$ is the state space, the drift term $f:\mathbb{R}^d\to\mathbb{R}^d$ and the diffusion term $g:\mathbb{R}^d\to\mathbb{R}^{d\times r}$ are Borel-measurable functions, and $B_t$ is a standard {$r${-dimensional} ({$r$-D})} Brownian motion defined on probability space $(\Omega,\mathcal{F},\{\mathcal{F}_t\}_{t\ge0},\mathbb{P})$ with a filtration $\{\mathcal{F}_t\}_{t\ge0}$ satisfying the regular conditions. The state $\vx^\ast$ with $f(\vx^\ast)=\vzero$ and $g(\vx^\ast)=\vzero$ is called equilibrium or zero solution.      

\paragraph{Notations.} Denote by $\Vert \cdot \Vert$ the Euclidean norm for a vector. Denote by $\Vert\cdot\Vert_{C(\mathbb{R}^d)}$ the maximum norm for the continuous function in $C(\mathbb{R}^d)$. $a\ll b$ means that the number $a$ is much smaller than $b$. For a function $V\in C^2(\mathbb{R}^d)$, denote by $\nabla V$ and $\mathcal{H}V$ the gradient and the Hessian matrix, respectively. 

% Denote by $\text{Lip}(\mathcal{X})$ the Lipschitz continuous space on $\mathcal{X}$.

%For $A=(a_{ij})$, a matrix of dimension $d\times r$, denote by $\|A\|^2_{\rm F}=\sum_{i=1}^{d}\sum_{j=1}^{r}a_{ij}^2$ the Frobenius norm.
\paragraph{Problem Statement} We assume that the zero solution of the Eq.~\eqref{SDE0}
is unstable, i.e. $\lim_{t\to\infty}\vx(t)\neq\vx^\ast$ on some set of positive measures. We aim to stabilize the zero solution using the neural controller with rigorous stability and safety guarantees that can be trained quickly. Specifically, we leverage the NNs to design an appropriate controller $\vu(\vx)$ with $\vu(\vx^\ast)=\vzero$ such that the controlled system 
\begin{equation}\label{SDE1}
	\mathrm{d} \vx(t)=f_\vu(\vx(t)) \mathrm{d} t+g(\vx(t))\mathrm{d} B_t
\end{equation}   
is steered to the zero solution with the whole controlled trajectory staying in the safety region, where $f_\vu\triangleq f+\vu$. Additionally, the neural controller should be optimized quickly in the training stage as well as in the testing stage such that it can be applied in the actual scenarios.

\paragraph{Basic Assumptions.} To guarantee the existence of unique solutions, it is sufficient for Eq.~\eqref{SDE1}
to own an initial condition and locally Lipschitz drift and diffusion terms. Therefore, we require the neural controllers to be locally Lipschitz continuous. This requirement is not overly onerous since the NNs can be designed by the compositions of Lipschitz continuous functions such as ReLU and affine functions. Without loss of generality, we assume that $\vx^\ast=\vzero$.

To investigate the analytical properties of the controlled SDE~\eqref{SDE1}, we introduce the stochastic derivative operator as follows.
\begin{definition}(\textbf{Derivative Operator})\label{derivative}
	Define the differential operator $\mathcal{L}_{\vu}$ associated with Eq.~(\ref{SDE1}) by
	\begin{equation*}
		\begin{aligned}
					\resizebox{1.0\linewidth}{!}{$
		\mathcal{L}_{\vu} \triangleq  \sum\limits_{i=1}^{d}(f_{\vu}(\vx))_i\dfrac{\partial}{\partial x_i}
		+\dfrac{1}{2}\sum\limits_{i,j=1}^{d}[g(\vx)g^{\top}(\vx)]_{ij}\dfrac{\partial^2}{\partial x_i\partial x_j}.
						$}
		\end{aligned}
	\end{equation*}
The subscript of $\mathcal{L}_{\vu}$ represents the dependence on $\vu$.
\end{definition}

According to the above definition of the derivative operator, an operation of $\mathcal{L}_{\vu}$ on
the function $V\in C^{2}(\mathbb{R}^d;\mathbb{R})$ yields:
\begin{equation*}\label{lie}
	\begin{aligned}
			\mathcal{L}V(\vx)=\nabla V({\vx})^\top f_\vu(\vx)
	+\dfrac{1}{2}\mathrm{Tr}\left[g^\top(\vx)\mathcal{H}V({\vx})g(\vx)\right].
	\end{aligned}
\end{equation*}
Next, we provide stability and safety analyses for Eq.~\eqref{SDE1} based on this derivative operator.

\subsection{Exponential Stability}\label{pre_stability}
In stochastic stability theory, an exponentially stable dynamical system implies that all solutions in some region around an equilibrium exponentially converge to this equilibrium almost surely (a.s.). Generally, stochastic stability theory investigates the exponential stability for the dynamics by studying the behavior of a potential function $V$ along the trajectories. Suppose the equilibrium is the minimum point of some potential function $V$. The convergence of the solution $\vx(t)$ is equivalent to the descent of the potential energy $V(\vx(t))$. In order to guarantee the exponential convergence, we can impose special structures for $V$ such that the potential field $\mathcal{L}_{\vu}V$ around the minimum is sinkable enough to attract the states from the high potential levels, as described in the following theorem.

\begin{thm}\label{thm1}
	\cite{mao2007stochastic}
	For the controlled SDE~\eqref{SDE1}, suppose that there exists a potential function $V\in C^{2}(\mathcal{X};\mathbb{R}_{+})$
	with $V(\vzero)=0$, constants $p>0$, $\varepsilon>0$, $c\in\mathbb{R}$ such that 
	\begin{equation}\label{eq in thm1}
		\begin{aligned}
			&		\varepsilon\|\vx\|^p\le V({\vx}),\\
			& \mathcal{L}_{\vu}V({\vx})\le cV({\vx}),
		\end{aligned}
	\end{equation}
	for all $\bm{x}\in\mathcal{X}$.
	Then,
	\begin{equation*}\label{exponentially stable}
		\limsup_{t\to\infty}\dfrac{1}{t}\log\|\vx(t;t_0,{\vx}_0)\|\le \dfrac{c}{p}~~a.s.,
	\end{equation*}
	for all solution $\vx(t;t_0,\vx_0)$ initiated from $(t_0,\vx_0)$. In particular, if $c<0$, the zero solution of Eq.~\eqref{SDE1}
	is exponentially stable almost surely.
\end{thm}

\paragraph{Learning Exponentially Stable Controller.} In order to learn the exponentially stable controller for SDEs,  Zhang et al. \citep{zhang2022neural} design a heuristic neural controller framework that integrates the conditions in Eq.~\eqref{eq in thm1} to the loss function. However, they can only assure the expected conditions are satisfied on finite training data instead of the whole state space. Hence, it lacks a rigorous stability guarantee and cannot be applied in stability-critical scenarios. In this paper, we circumvent this drawback for heuristic neural controllers through a projection operator.

\paragraph{Difference From the Lyapunov Theory}
One may naturally regard the potential function $V$ and the Theorem~\ref{thm1}, respectively, as the classic Lyapunov function and the Lyapunov stability theory for ODEs. The main difference is that the Lyapunov stability guarantees the decrease of the potential function, but the stochastic theory can only guarantee that the potential function eventually converges to a minimum. In other words, trajectories of exponentially stable SDEs may reach undesired states due to the stochasticity, which is not safe as we elaborate in the following subsection.

\subsection{Safe Controller Synthesis}\label{pre safe controller}
Safety-critical systems in application domains, such as automatic drive and medicine, require their state trajectories to stay in the specified safe region in order to prevent loss of life and economic harm. Here, we formalize the concept of safety in the following definition.
\begin{definition}(\textbf{Safe Controller})
The controller in the controlled system~\eqref{SDE1} is a safe controller if every solution $\vx(t)$ initiated from $\vx(0)\in\rm{int}(\mathcal{C})$ satisfying:
\begin{equation*}
	\forall t\ge0:~\vx(t)\in\mathcal{C}~a.s.,
\end{equation*}
 where $\mathcal{C}$ is
a specified safe region induced by a
locally Lipschitz function $h$ as $\mathcal{C}=\{\vx:h(\vx)\ge0\}$ and $\vzero\in\mathcal{C}$.
\end{definition}
Recent works on the safe controller synthesis typically treat $\partial\mathcal{C}=\{\vx:h(\vx)=0\}$ as a barrier, and require the controlled dynamics to satisfy additional conditions such that the trajectories never cross the barrier~\citep{clark2019control}. In addition, existing works design complicated barrier functions, and/or  the controller through the online quadratic program (QP) method, which requires a high computational cost~\cite{fan2020bayesian,sarkar2020high}. In Section \ref{section safety}, we propose a theorem for safety guarantee with a simpler barrier function and design an offline safe controller based on this result.

\paragraph{Connection to Safe RL} 
The flow maps $\{S_t\}_{t\ge0}$ of the Eq.~\eqref{SDE1} can be modeled as:
\begin{equation}\label{flow}
S_t(\vx_0)\triangleq\vx_0+\int_0^tf_{\vu}(\vx(s))\mathrm{d}s+\int_0^tg(\vx(s))\mathrm{d}B_s,
\end{equation}
which induce a class of Markovian decision processes (MDPs)~\citep{filar2012competitive}. Hence, the safe controller synthesis with a fixed time point $t$ can be regarded as the safe reinforcement learning (RL) problems~\citep{garcia2015comprehensive}. The major gap for applying RL methods to our problem is that, generally, the analytical expression of $S_t$ is intractable given Eq.~\eqref{SDE1}. And the unbounded diffusion term $g$ further incurs difficulty for handling safety issues in RL. Besides,
existing RL methods consider the safety in the meaning of average over trajectory or risk probability~\citep{chow2017risk,sootla2022saute}, rather than satisfying the safety constraint almost surely (or with the probability one) as proposed in this paper.

\section{Neural Controller with Exponential Stability Guarantee}\label{section stability}
%We now present LyaNet, our Lyapunov framework for training ODEs of the form specified by Equations (1) to (3). As
%alluded to in Section 2.3, our goal is to find parameters θ
%of the ODE to satisfy the Lyapunov exponential stability
%condition in Theorem 1 with respect to a potential function
%V . We develop the formulation in two steps:
%1. Section 3.1: For a given supervised loss, we define an
%appropriate potential function V .
%2. Section 3.2: For that V , we define the Lyapunov loss
%which captures the degree of violation from satisfying
%the contraction condition in Equation (10) that implies
%exponential stability.
%Theoretically, we show that optimizing the Lyapunov loss
%implies the learned ODE exponentially stabilizes to predictions with minimal supervised loss (Theorem 2), which in
%turn implies a novel adversarial robustness guarantee (Theorem 3).6
%In other words, if we find a θ ∈ Θ that achieves
%zero Lyapunov loss, then Equation (10) will be satisfied. We
%present practical learning algorithms in Section 4.

We leverage the NNs to learn the neural controller with the exponential stability guarantee, i.e. satisfying the exponential stability condition as specified in Theorem~\ref{thm1}. 
The specific process is divided into two steps: 1) we propose a heuristic neural framework for learning the neural controller according to an optimization problem based on Theorem~\ref{thm1}; 2) A projection operator is developed to pull the learned controller into the stable domain, and we provide an operator with a closed form to approximate this projection operator in Theorem~\ref{thm2}. Furthermore, we analyze the rationality of the conditions in Theorem~\ref{thm2} based on the safety of the controller.

\subsection{Heuristic Neural Controller}\label{step1 stability}
A consensus in the control field is that the controllers should be as simple and small as possible~\citep{lewis2012optimal}. Hence, we use a class of Lipschitz functions on $\mathcal{X}$, denoted as $\text{Lip}(\mathcal{X})$, to represent the space of candidate controllers. Incorporating with the Theorem~\ref{thm1}, we can obtain the following optimal problem for learning the stabilization controller:
\begin{equation*}
	\begin{aligned}
		\arg&\min_{\vu\in\text{Lip}(\mathcal{X})} \frac{1}{2}\Vert\vu(\vx)^\top R\vu(\vx)\Vert_{C(\mathbb{R}^d)}\\
		\text{s.t.}~&\varepsilon\Vert\vx\Vert^p-V(\vx)\le0,\\
		 &\mathcal{L}_{\vu}V(\vx)-cV(\vx)\le0,
	\end{aligned}
\end{equation*}
where $\varepsilon>0,~c<0,~p>0$ are predefined hyper-parameters, and semi-positive-definite matrix $R$ measures the importance of different control components.

\paragraph{Neural Network Framework.} In order to find the closed-form controller $\vu$, we transform the above optimization problem into an NN framework. Specifically, we parameterize the controller by an NN, i.e.,  $\vu=\vu_{\vtheta}(\vx)$ with $\vu_{\vtheta}(\vzero)=\vzero$, where $\vtheta$ is trainable parameter vector. Additionally, we constrain the Lipschitz constant of the neural controller by utilizing the spectral normalization method~\citep{yoshida2017spectral,miyato2018spectral}. The other Lipschitz neural networks can be used as well \citep{virmaux2018lipschitz,liu2022learning}. Then, we construct the neural potential function $V_{\vtheta}$ by using the convex functions~\citep{icnn} such that it satisfies the positive definite condition, i.e., $V_\vtheta(\vx)\ge0$ and $V_\vtheta(\vzero)=0$. Second-order differentiable activation functions are used to guarantee $V_{\vtheta}\in C^{2}(\mathbb{R};\mathbb{R}_{+})$. We add $L^p$ regularization term $\varepsilon\Vert\vx\Vert^p$ with $\varepsilon\ll1$ to $V_\vtheta$ such that  $V_\vtheta(\vx)\ge \varepsilon\Vert\vx\Vert^p$. Then, we define the supervised stabilization loss function as follows.
\begin{definition}(\textbf{Stabilization Loss}) Consider a candidate potential function $V_\vtheta$ and a controller $\vu_\vtheta$ for the
	controlled system~\eqref{SDE1}. The exponential stabilization loss is defined as
\begin{equation}\label{loss1}
	\begin{aligned}
		L_{es}(\vtheta)\triangleq\frac{1}{N}\sum_{i=1}^N\big[\vu_\vtheta&(\vx_i)^\top R\vu_\vtheta(\vx_i)\\
		+\lambda_1\max\big(0,\mathcal{L}_{\vu_{\vtheta}}&V_\vtheta(\vx_i)-cV_\vtheta(\vx_i)\big)\big],
	\end{aligned}
\end{equation}
where $\lambda_1>0$ is a hyper-parameter representing the weight factor and $\{(\vx_i)\}_{i=1}^N$ is the dataset.

\end{definition}

\paragraph{Weaknesses of The Heuristic Framework.}  
Two aspects could be improved in the current learning framework. First, the NNs trained on finite samples cannot guarantee the condition in the loss function is also satisfied in the infinite set $\mathcal{X}$. Second, the computational complexity of computing the Hessian matrix $\mathcal{H}V$ is $\mathcal{O}(d^2)$, a major bottleneck for applications in high-dimensional systems. We introduce the projection operation in the following subsection to overcome the first weakness. We defer the solution to the second weakness in Section~\ref{acceleration section}.

\subsection{Projection Operator for Exponential Stability}\label{step2 stability}
We define the space of the exponentially stable controllers as follows:
\begin{equation*}
	\mathcal{U}_{es}(V,\mathcal{X})\triangleq\{\vu\in \text{Lip}(\mathcal{X}):\mathcal{L}_{\vu}V-cV\le0\}.
\end{equation*}
To force the learned controller in $\mathcal{U}_{es}(V,\mathcal{X})$, a common technique is to find the $L^2$ projection of the controller which is defined as follows.
\begin{definition}(\textbf{Projection Operator})
	Denote by $\pi(\vu,\mathcal{U})$ the projection from $\vu$ onto the target space $\mathcal{U}$ as 
	\begin{equation*}
		\pi(\vu,\mathcal{U})\in\arg\min_{\tilde{u}\in\mathcal{U}}\frac{1}{2}\Vert\vu(\vx)-\tilde{\vu}(\vx)\Vert^2_{C(\mathbb{R}^d)}.
	\end{equation*}
\end{definition}
In General, it is difficult to obtain a closed form of the above projection operator. Recently, the QP method has been used to solve the following simplified projection for each data $\vx$ as an alternative~\citep{chow2019lyapunov},
\begin{equation*}
	\begin{aligned}
		&\tilde{\pi}(\vu)(\vx)\in\arg\min\frac{1}{2}\Vert\vu(\vx)-\tilde{\vu}(\vx)\Vert^2,\\
		&~~~~~~~~~~~~~~~~\text{s.t.}~\mathcal{L}_{\tilde{\vu}}V(\vx)-cV(\vx)\le0.
	\end{aligned}
\end{equation*} 
However, this simplified operation $\tilde{\pi}$ cannot assure the continuity of the projection and causes additional computational costs. In the following theorem, we approximate the projection operator $\pi$ with a closed form.

\begin{figure}[t]
% 	\vskip -0.1 in
	\centering
	\includegraphics[width=8.2cm]{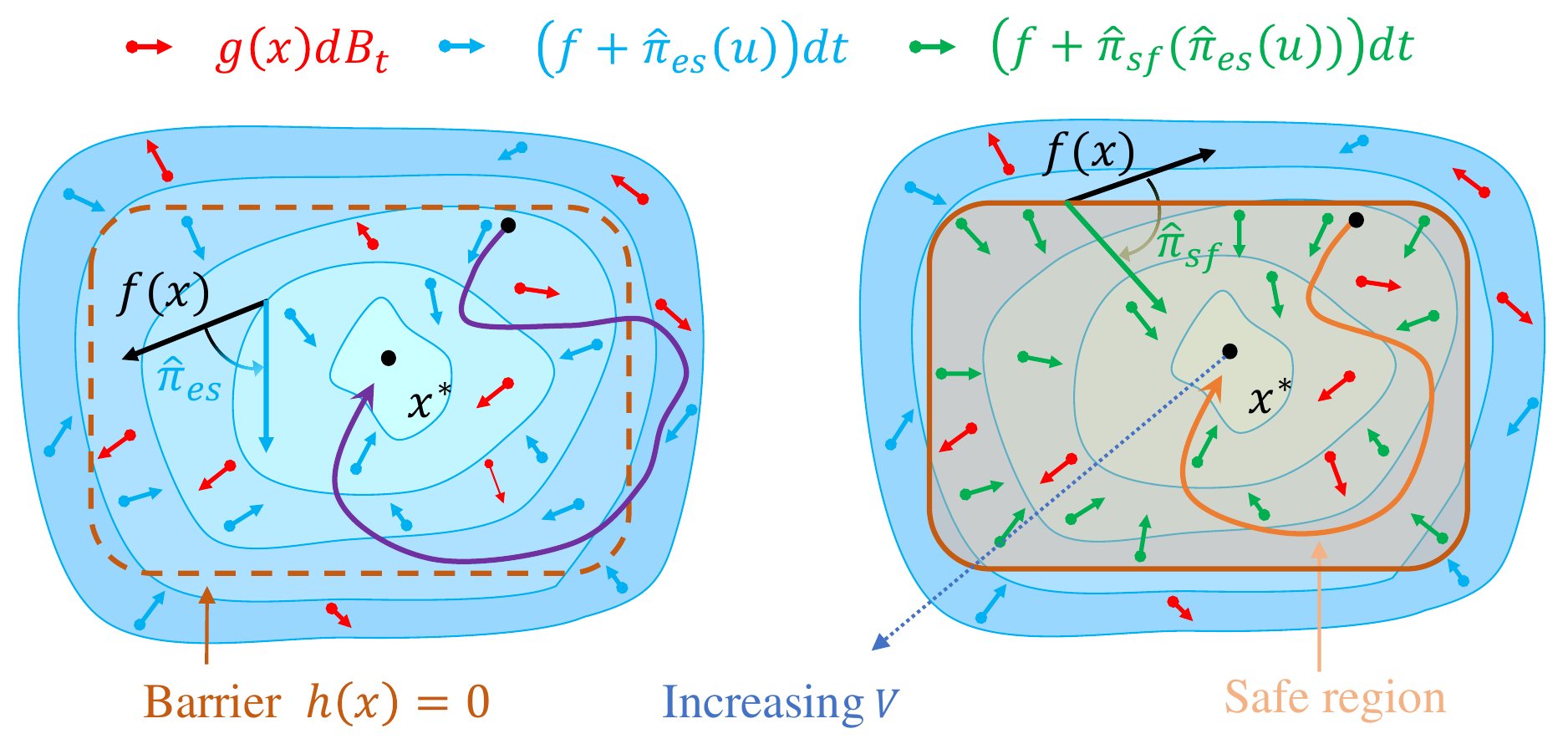}
	\caption{ Illustration of the stable and safe controlled  dynamics. \textbf{Left:} Operator $\hat{\pi}_{es}$ projects candidate controller $\vu$ to the exponentially stable domain, i.e. $\mathcal{L}_{\vu}V\le-cV$, s.t. any trajectory will exponentially converge to the equilibrium point $\vx^{\ast}$. The trajectory cannot stay in the safe region (green shaded area) due to the stochasticity. \textbf{Right:} Operator $\hat{\pi}_{sf}$ project the controller to the safe domain, i.e. $\mathcal{L}_{\vu}h\ge-\alpha(h)$, to restrict the trajectory (orange line) in the safe region.
	}
	\label{overflow}
% 	\vskip -0.25 in
\end{figure}

\begin{thm}\label{thm2}(\textbf{Approximate Projection Operator})
	For a candidate controller $\vu$ and the target space $\mathcal{U}_{es}(V,\mathcal{X})$, an approximate projection operator is defined as follows,
	\begin{equation}\label{st proj}
		\hat{\pi}_{es}(\vu,\mathcal{U}_{es}(V,\mathcal{X}))\triangleq\vu-\dfrac{\max(0,\mathcal{L}_{\vu}V-cV)}{\Vert \nabla V\Vert^2}\cdot\nabla V,
	\end{equation}
 then, the operator $\hat{\pi}_{es}$ satisfies  $\hat{\pi}_{es}(\vu,\mathcal{U}_{es}(V,\mathcal{X}))\in\mathcal{U}_{es}(V,\mathcal{X})$ if and only if the state space $\mathcal{X}$ of the controlled systems~\eqref{SDE1} is bounded.
\end{thm}
The proof is shown in the Appendix~\ref{proof thm2}. We can directly apply this approximate projection operator to the learned controller to provide the exponential stability guarantee when $\mathcal{X}$ is bounded. However, this bounded condition may not hold on a general exponentially stable system, the following example illustrates this point.
\begin{ex}\label{ex1}(\textbf{Exponential Martingale})
	Consider a 1-D SDE described by $
		\mathrm{d}x=ax\mathrm{d}t+bx\mathrm{d}B_t,$
where $a<0,~|a|>\frac{1}{2}b^2$ and $x_0>0$, the zero solution is exponentially stable, but for any $M>0$, $\mathbb{P}(x(t)>M)>0$, i.e. $\mathcal{X}=\mathbb{R}$.
\end{ex}
See the poof in the Appendix~\ref{proof ex1}. This example implies that the boundedness of state space $\mathcal{X}$ cannot be induced under the exponential stability. As introduced in Section~\ref{pre safe controller}, if the controller is safe in the bounded safe region $C$, then $\mathcal{X}$ is bounded for any $\vx_0\in\mathcal{C}$. %This inspires us that to bound the $\mathcal{X}$,
We, therefore, take safety into consideration to bound the $\mathcal{X}$. Furthermore, this example shows that a stable controller may not be safe.

\section{Neural Controller with Safety Guarantee}\label{section safety}
In this section, we present an approach for safe controller synthesis. Similar to the analysis in Section~\ref{section stability}, the safe controller's construction can be divided into two steps as well: 1) training a heuristic neural safe controller (Section~\ref{step1 safe}), and 2) projecting the learned controller to obtain the safety guarantee in the post-training stage (Section~\ref{step2 safe}).  In particular, we propose a new sufficient condition for the safe controller based on the barrier function theory~\citep{ames2016control}. Then one can design a neural safe controller according to this result. In addition, we show that the intersection of the space of the safe controller and $\mathcal{U}_{es}$ is not empty, which assures the existence of a controller with both safety and stability guarantees.

\subsection{Zeroing Barrier Function for Safety}\label{step1 safe}
As specified in Section~\ref{pre safe controller}, the primary task for safe controller synthesis is to impose additional structures for the controlled system~\eqref{SDE1} to restrict the trajectory, never approaching the barrier $\partial \mathcal{C}$. Here, we present a concise structure for safety based on the zeroing barrier function.
\begin{definition}\citep{ames2016control}
	A continuous function $\alpha: (0, a)\to(-\infty,\infty)$
	is said to belong to class-$\mathcal{K}$ for some $a > 0$ if it
	is strictly increasing and $\alpha(0) = 0$.
\end{definition}
\begin{thm}\label{thm3}
	 For the SDE (\ref{SDE1}), if there exists a class-$\mathcal{K}$ function $\alpha(x)$ such that
	\begin{equation}\label{eq for safety}
		\mathcal{L}_{\vu}h(\vx)\ge-\alpha(h(\vx))
	\end{equation}
	for $\vx\in\mathcal{C}$. Then, the controller $\vu$ is a safe controller, and $h(\vx)$ is called a zeroing barrier function (ZBF).
\end{thm}
The proof is provided in Appendix~\ref{proof thm3}. An intuition behind this theorem is that 
 $\alpha(h)$ can be regarded as a potential function. When the trajectory $\vx(t)$ approaches the barrier $\partial\mathcal{C}$, the potential $\alpha(h(\vx(t)))$ decreases to the minimum $0$ so that $\mathcal{L}_{\vu}h\ge0$, which in turn lifts the potential to pull the trajectory back
to the interior of $\mathcal{C}$.
\paragraph{Neural Safe Controller.} Given the ZBF $h(\vx)$, we can also formulate the safety problem as an optimization problem, 
\begin{equation*}
	\begin{aligned}
	\arg&\min_{\vu\in\text{Lip}(\mathcal{C}),\alpha\in\mathcal{K}} \frac{1}{2}\Vert\vu(\vx)^\top R\vu(\vx)\Vert_{C(\mathbb{R}^d\times \mathbb{R}_{+})}\\
	\text{s.t.}~&-\mathcal{L}_{\vu}h(\vx)-\alpha(h(\vx))\le0.
\end{aligned}
\end{equation*}
To convert this problem into a learning-based framework, we adopt the same neural controller $\vu_{\vtheta}$ as in Section~\ref{step1 stability}. We apply the monotonic NNs to construct the candidate extended class-$\mathcal{K}$ function as $\alpha_{\vtheta}(x)=\int_0^xq_{\vtheta}(s)\mathrm{d}s$, where $q_{\vtheta}(\cdot)$, an NN, is definitely positive~\citep{umnn}. Then we train the monotonic NNs with the safe loss as follows.
\begin{definition}(\textbf{Safety Loss}) For a candidate class-$\mathcal{K}$ function $\alpha_\vtheta$ and a controller $\vu_\vtheta$, the safety stabilization loss is defined as
\begin{equation}\label{loss2}
	\begin{aligned}
	L_{sf}(\vtheta)\triangleq&\frac{1}{N}\sum_{i=1}^N\big[\vu_\vtheta(\vx_i)^\top R\vu_\vtheta(\vx_i)\\
	+\lambda_2&\max\big(0,-\mathcal{L}_{\vu_{\vtheta}}h(\vx_i)-\alpha_\vtheta(h(\vx_i))\big)\big],
	\end{aligned}
\end{equation}
where $\lambda_2>0$ is a tunable weight factor and $\{(\vx_i)\}_{i=1}^N$ is the dataset.

\end{definition} 

\subsection{Projection Operator for Safety}\label{step2 safe}
Notably, the learned neural safe controller may not satisfy the expected condition
in Eq.~\eqref{eq for safety}. However, we can project it to the safe controller's space $\mathcal{U}_{sf}(\alpha,\mathcal{C})\triangleq\{\vu\in\mathcal{U}(\mathcal{C}):-\mathcal{L}_{\vu}h-\alpha(h)\le0\}$ 
as $\pi(\vu_\vtheta,\mathcal{U}_{sf}(\alpha_\vtheta,\mathcal{C}))$.  To efficiently implement the projection, we construct an approximate projection operator with a closed form again.
\begin{thm}\label{thm4}
	For the bounded safe region $\mathcal{C}$, the approximate projection operator defined as
	\begin{equation}\label{sa proj}
	\hat{\pi}_{sf}(\vu,\mathcal{U}_{sf}(\alpha,\mathcal{C}))\triangleq\vu+\dfrac{\max(0,-\mathcal{L}_{\vu}h-\alpha (h))}{\Vert \nabla h\Vert^2}\cdot\nabla h,
	\end{equation}
	 satisfies  $	\hat{\pi}_{sf}(\vu,\mathcal{U}_{sf}(\alpha,\mathcal{C}))\in\mathcal{U}_{sf}(\alpha,\mathcal{C})$.
\end{thm}
\begin{remark}
	To maintain $\hat{\pi}_{sf}(\vu,\mathcal{U}_{sf}(\alpha,\mathcal{C}))(\vzero)=\vzero$, we further require $\nabla h(\vzero)=\vzero$, this condition is easy to satisfy since one can modify the ZBF around the equilibrium without change the safe region $\mathcal{C}$.
\end{remark}
The proof is shown in Appendix~\ref{proof thm4}. Based on this approximate projection operator in the post-training stage, one can transform the learned controller into a safe controller $\hat{\pi}_{sf}(\vu_\vtheta,\mathcal{U}_{sf}(\alpha_\vtheta,\mathcal{C}))$. 

\paragraph{Controller with Both Safety and Exponential Stability.}  Actually, we can jointly train the $\vu_\vtheta, V_{\vtheta},\alpha_\vtheta$  with the summed loss $L_{es}(\vtheta)+L_{sf}(\vtheta)$ to learn a safe and exponentially stable controller with the finite sample size. A natural question is whether we can composite the approximate projection operator $\hat{\pi}_{es}$ and $\hat{\pi}_{sf}$ together to get the controller with rigorous safety and stability guarantees or not. To answer this question, we study the relationship between the safe controller space $\mathcal{U}_{sf}(\cdot,\mathcal{C})$ and the exponentially stable controller space $\mathcal{U}_{es}(\cdot,\mathcal{C})$.

\begin{thm}\label{thm5}(\textbf{Existence of the Safe and Exponentially Stable Controller})
	For any ZBF $h(\vx)$ with a bounded safe region $\mathcal{C}$ and a class-$\mathcal{K}$ function $\alpha$, there exists a potential function $V$, s.t. $\mathcal{U}_{es}(V,\mathcal{C})\cap\mathcal{U}_{sf}(\alpha,\mathcal{C})\neq\emptyset$.
\end{thm}

\paragraph{Proof Sketch.} For the special case that $h(\vzero)=\max_{\vx\in\mathcal{C}}h(\vx)$, we can construct the potential function as $V(\vx)=h(\vzero)-h(\vx)$ s.t. $V\ge k\|\vx\|$ for some $k>0$. Then it can be shown that there exists $\vu\in\mathcal{U}_{sf}(\alpha,\mathcal{C})$ s.t. $\mathcal{L}_{\vu}V\le -kV$. If $h(\vzero)\neq\max_{\vx\in\mathcal{C}}h(\vx)$, we can modify the ZBF as $\tilde{h}=h+\lambda$ s.t. $\tilde{h}(\vzero)=\max_{\vx\in\mathcal{C}}\tilde{h}(\vx)$ without changing the safe region, where $\lambda$ is a smooth approximation of the Dirac function $\delta_{\vzero}$. Similar to the special case, we can derive the remaining parts of the theorem. The detail of the proof is shown in Appendix~\ref{proof thm5}.
{\color{black}
\begin{remark}\label{rem1}
    The bounded safe region condition is Theorem~\ref{thm5} can be relaxed to unbounded cube $[a_1,b_1]\times\cdots\times[a_d,b_d]$, $-\infty\le a_i<b_i\le\infty$, where at least one interval $[a_i,b_i]$ is bounded. The proof is provided in Appendix~\ref{proof thm5}.
\end{remark}}
According to Theorem~\ref{thm5}, we can first apply operator $\hat{\pi}_{sf}$ to the learned controller $\vu_\vtheta$ 
to obtain the safety guarantee and restrict the state space $\mathcal{X}\subset\mathcal{C}$. Since $\mathcal{X}$ is bounded now, then, we can apply operator $\hat{\pi}_{es}$ to the safe controller to obtain the exponential stability guarantee according to Theorem~\ref{thm2}.

\section{Accelerate the Training and the Projection}\label{acceleration section}
In general, computing $\mathcal{L}_{\vu}V$ costs $\mathcal{O}(d^2)$. The term $\mathrm{Tr}[g^\top\mathcal{H}Vg]$ requires computing a Hessian matrix of $V$. The same is valid for computing $\mathcal{L}_{\vu}h$. Here, two tricks can be used to achieve low computational costs. First, due to the commutativity of the trace operator, we can get the following unbiased estimate:
\begin{equation*}
	\mathrm{Tr}[g^\top\mathcal{H}Vg]=\mathbb{E}[\xi^\top\mathcal{H}Vgg^\top\xi]=	\mathbb{E}[(\nabla(\xi^\top\nabla V))^\top gg^\top\xi],
\end{equation*}
where $\xi$ is a $d$-D noise vector with mean zero and variance $I$. The Monte Carlo estimator for the above expectation is known as the
Hutchinson’s trace estimator~\citep{hutchinson1989stochastic}.  Second, if $g$ is a vector independent of $\vx$, computing the gradient of the vector-Jacobian product $g^\top\nabla V$ costs $\mathcal{O}(d)$, and the remaining part is simply a vector-Jacobian product $\mathrm{Tr}[g^\top\mathcal{H}Vg]=g^\top\nabla(g^\top\nabla V)$. These two tricks can be used according to the shape of the diffusion term $g(\vx,t)\in\mathbb{R}^{d\times r}$, i.e.,  matrix ($r>1$) or vector ($r=1$).

\paragraph{Case 1: $r>1$.} In this case, $g$ is a matrix, and we use the former trick to accelerate the training stage.
Assuming the cost of evaluating $\mathrm{Tr}[g^\top\mathcal{H}Vg]$ is on the order of $\mathcal{O}(Nd^2)$ where $d$ is the dimensionality of the
data and $N$ is the size of the dataset, then the cost of computing the following Hutchinson’s trace estimator via automatic
differentiation, which only involve vector-Jacobian products, is $\mathcal{O}(NMd)$, 
\begin{equation}\label{hutch1}
	\mathrm{Tr}[g^\top\mathcal{H}Vg]\approx\frac{1}{M}\sum_{i=1}^{M}(\nabla(\xi_i^\top\nabla V))^\top gg^\top\xi_i,
\end{equation}
 where $M$ is the number of samples in the Monte Carlo estimation with $M\ll d$. In practice,  $M$ can be tuned to trade off the variance and the computational cost. As suggested in~\citep{grathwohl2018ffjord,song2020sliced}, $M = 1$ is already a good choice, which worked well in our experiments. Notice that although the approximate projection operator applied to the learned controller in the projection stage also involves computing the $\mathcal{L}_{u}V$, we cannot use this trick here because it does not equal the real trace.
 
\paragraph{Case 2: $r=1$.}
Since $g$ is a vector, we can use the latter trick by removing the dependency of $g$ on data $\vx$. Specifically, after getting the $g(\vx)$ on data $\vx$, we apply the $\texttt{detach}$ operator in Pytorch~\citep{paszke2019pytorch} to $g$ to prevent calculating the gradients w.r.t $\vx$. Then the computational cost reduces from $\mathcal{O}(Nd^2)$ to $\mathcal{O}(Nd)$. This trick is shown as follows:
\begin{equation}\label{hutch2}
	\mathrm{Tr}[g^\top\mathcal{H}Vg]=g^\top\nabla((g.\texttt{detach()})^\top\nabla V).
\end{equation}
We note that the result in Eq.~\eqref{hutch2} is an identity equation. Hence, it is an ideal choice when $r=1$. Since the derivative operator $\mathcal{L}_{\vu}$ is calculated in both the training and post-training projection stages, and these two cases always cover at least one stage, we can accelerate the whole framework of the proposed safe and exponentially stable controller. We summarized the complete framework of the proposed FESSNC in Algorithm~\ref{algo1}.

\begin{algorithm}[t]
	\centering
	\caption{ \footnotesize{FESSNC} } \label{algo1}
	\begin{algorithmic}
		\STATE {\bfseries Input:} Dynamics $f,g$, ZBF $h$, safe region $\mathcal{C}$, max iterations $m$, learning rate $\gamma$, initial parameters $\vtheta$, weight factors $\lambda_{1,2}$, noise vector samples $\{\xi_i\}_{i=1}^M$.
		\STATE {\bfseries Output:}   Controller $\vu_\vtheta$, potential function $V_\vtheta$, class $\mathcal{K}$ function $\alpha_\vtheta$.
		\FOR{$r=1:m$}
		\STATE $\{\vx_i\}_{i=1}^N\sim\mathcal{C}$  \hfill {$\triangleright$ Sample training data}
		\STATE $L(\vtheta)=L_{es}(\vtheta)+L_{sf}(\vtheta)$   \hfill $\triangleright$ Compute loss \eqref{loss1}\eqref{loss2}\eqref{hutch1}
		\STATE $\vtheta\gets \vtheta - \gamma \cdot   \nabla_{\vtheta} L({\vtheta})$   \hfill {$\triangleright$ Update parameters}
		\ENDFOR
		\STATE {\bfseries Return:} $\hat{\pi}_{es}(\hat{\pi}_{sf}(\vu_\vtheta,\mathcal{U}_{sf}(\alpha_\vtheta,C)),\mathcal{U}_{es}(V_\vtheta,\mathcal{C}))$ 
		\STATE \hfill {$\triangleright$ Safety and exponential stability guarantee \eqref{st proj}\eqref{sa proj}\eqref{hutch2}}
	\end{algorithmic}
\end{algorithm}

\section{Experiments}\label{experiments}
% In this section, we demonstrate the power of the proposed FESSNC on several case studies as well as its application in synchronizing the coupled oscillators. 
{\color{black}In this section, we demonstrate the superiority of the FESSNC over existing methods using several case studies, then we validate the scalability of FESSNC to the high dimensional system with an application in synchronizing the coupled oscillators. More details of the experiments can be found in Appendix~\ref{appen_details}.}
\subsection{Case Studies}
 We first show that our approach is able to produce reliable control policies with stability and safety guarantees for each task. We then show the advantages of the FESSNC in computational costs compared with representative baselines. 

\begin{figure}[htp]
	\centering
	\includegraphics[width=6.5cm]{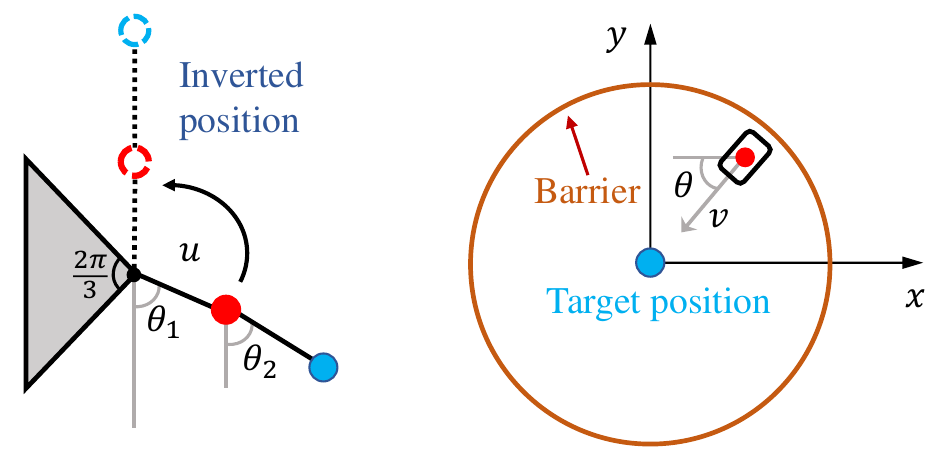}
	\caption{ Schematic diagrams of the double pendulum (\textbf{left}) and the kinetic bicycle (\textbf{right}) with state constraints.
	}
	\label{sketch}
	\vskip -0.15 in
\end{figure} 
\paragraph{Models}
We consider two classic control models, the fully actuated double pendulum swing-up and the fully actuated planar kinetic bicycle motion. The double
pendulum is a $4$-D system with $2$ angles and $2$ angular velocities, $(\theta_1,\dot{\theta}_1,\theta_2,\dot{\theta}_2)$~\citep{deisenroth2013gaussian}. We add multiplication noise force to the system such that the angular acceleration is perturbed related to angles.  The objective is to learn a
control strategy that steers the double pendulum to the inverted position within the safe region where the inner pendulum can only spin in  $[-\pi/6,7\pi/6]$, see Figure~\ref{sketch} (left). The kinematic bicycle model for car-like vehicles has $4$ variables $(x, y, \theta, v)$, where $(x,y)$ is the bicycle's position, $\theta$ is the direction of motion, and $v$ is the velocity~\citep{rajamani2011vehicle}. The goal is to steer the noise-perturbed bicycle to the target position without crossing the round wall with a radius of $2$, see Figure~\ref{sketch} (right).

\begin{figure}[b]
	\centering
	\includegraphics[width=8.2cm]{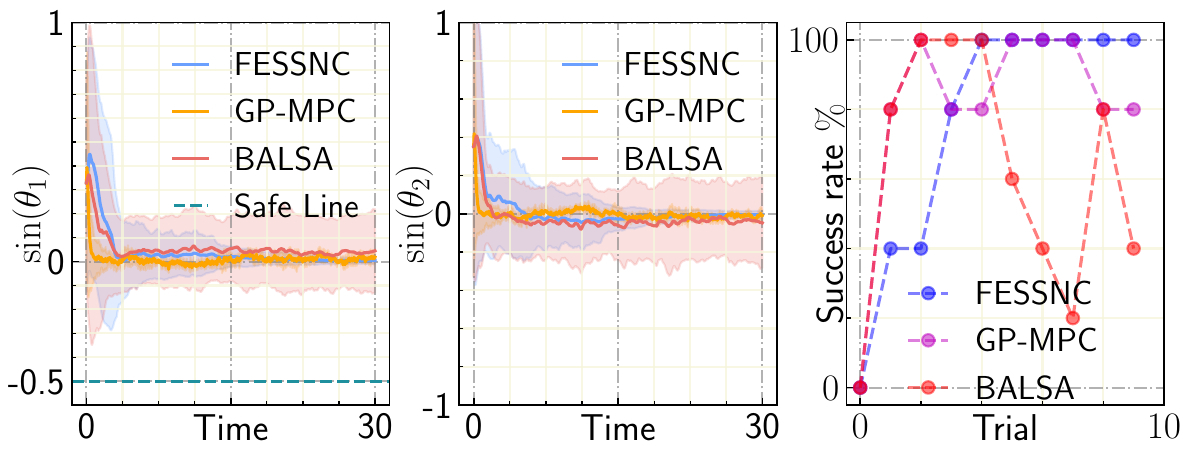}
	\caption{Results of the constraint double pendulum via different methods. (\textbf{left} resp., \textbf{center})) Sine values
		 along the trajectories of $\theta_1$ (resp., $\theta_2$). Safe region is $\sin(\theta_1) \le\frac{1}{2}$. (\textbf{right}) The success rate of the three methods is $3$ seconds per trial.
	We define 'success' if the pendulum tip's angle is closer than
	$\frac{\pi}{40}$ to the upright position for consecutive $3$ seconds.}
	\label{double pendulum}
\end{figure} 

\paragraph{Baselines}
To simulate the trajectories of the controlled system, we approximate the flow $S_t$ of the SDEs~\eqref{flow} with the Euler–Maruyama method~\citep{sarkka2019applied}. Then the obtained MDPs (or discrete dynamics) can be handled with model-based RL methods. We compare the proposed FESSNC with the GP-MPC algorithm~\citep{kamthe2018data} and the BALSA algorithm~\citep{fan2020bayesian}, respectively the most data-efficient MPC-based and QP-based RL algorithms to date.

\begin{figure}[t]
	\centering
	\includegraphics[width=8.2cm]{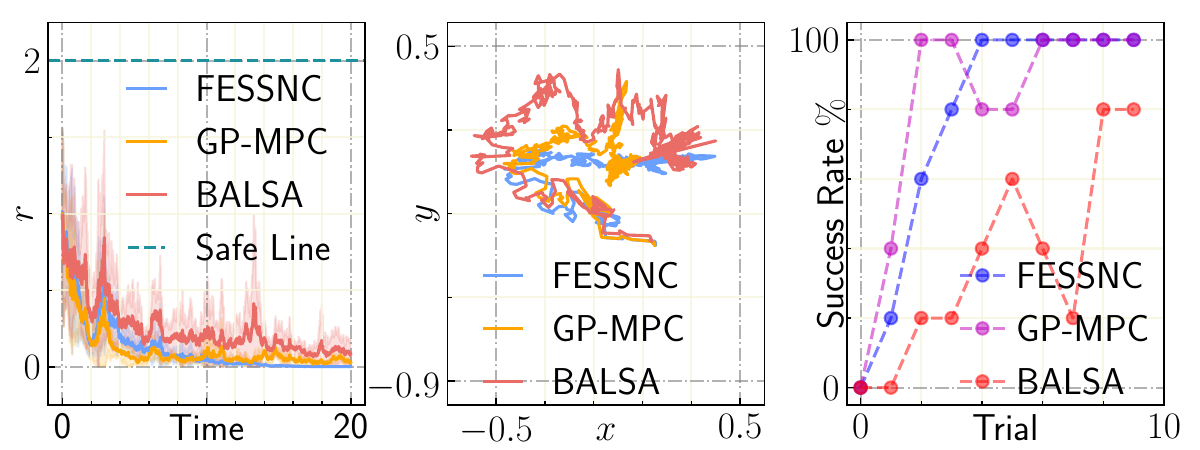}
	\caption{Results of the constraint kinematic bicycle model. (\textbf{left}) The distance $r$ between the car and the target position. The safe region is $x^2+y^2\le 4$. (\textbf{center}) The mean of motion trajectories. (\textbf{right}) The success rate of the three methods is $3$ seconds per trial. We define 'success' if the car is closer than
		$0.1$ to the target position for consecutive $2$ seconds.
	}
	\label{bicycle}
\end{figure}

We conduct $5$ independent experiments, and the results are shown in Figure~\ref{double pendulum}-\ref{bicycle}, and Table~\ref{table1}. In the figures, the solid lines are obtained by averaging
the sampled trajectories, while the shaded areas stand
for the variance regions. 

\paragraph{Performance.}  As shown in Figure~\ref{double pendulum}, all three methods satisfy the safety constraint for the double pendulum task. The FESSNC and GP-MPC can quickly steer the pendulum to the inverted position, but the GP-MPC cannot balance the pendulum abidingly, whereas the FESSNC does. The BALSA performs poorly because it requires the potential function $V$ and ZBF $h$ to be defined prior, and the predefined functions may induce disjoint safe controllers space and exponentially stable controllers space. The results in the kinetic bicycle task show the same situations as the double pendulum in Figure~\ref{bicycle}. The motion trajectories of the bicycle under the FESSNC are the most robust.

\paragraph{Computational Cost.} We compare the training time of one complete control process for the FESSNC and baselines. The results are shown in Table~\ref{table1}. We can see that the GP-MPC costs much more time than the others. This is because the GP-MPC method has to roll out the forward dynamics to obtain the predictive and adjoint states, and then find the optimal control through gradient descent. The BALSA is slightly quicker than FESSNC, but its performance is far inferior to the latter. Moreover, QP-based solver costs at least $\mathcal{O}(d^2)$ for $d$-D system~\citep{kouzoupis2018recent}, which makes it hard to be applied in the high-dimensional system, whereas our FESSNC does.

\begin{table}[htb]
	\vskip -0.1in
\caption{Average amount of time in seconds needed to complete a control process for each method.}\label{table1}
	\vskip 0.1in
\centering
\footnotesize{
\begin{tabular}{ccc}
		\toprule
		\toprule
		Method
		&{Double Pendulum}
		&{Kinetic Bicycle}
		\\ 
	\midrule 
		FESSNC   &$10.35$s  &$6.39$s \\
		GP-MPC   &$146.64$s   &$26.68$s \\
		BALSA   &$3.08$s    &$2.18$s \\
				\midrule
				\bottomrule
	\end{tabular}
 }
\end{table} 

\begin{figure}[htp]
	\centering
	\includegraphics[width=8.4cm]{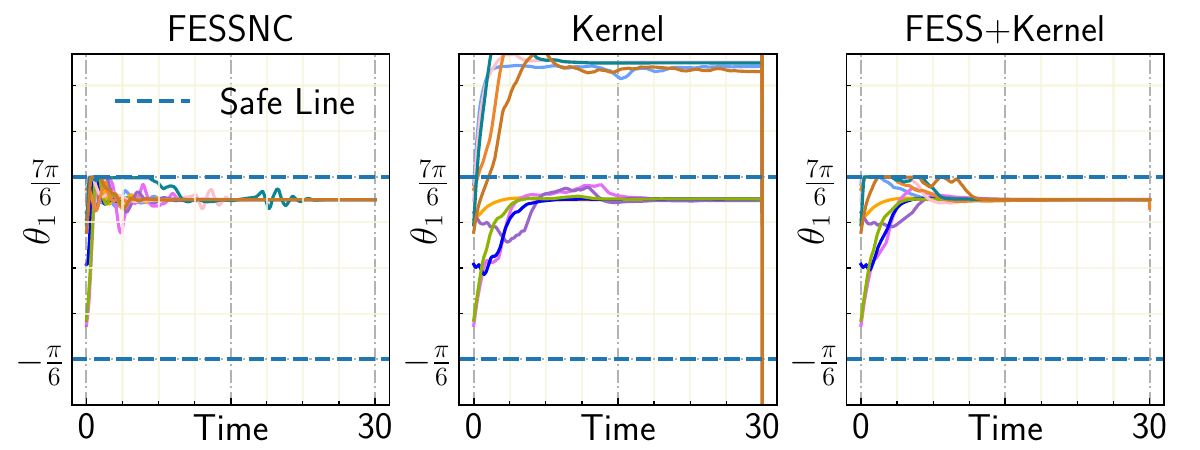}
	\caption{{\color{black}The angle of the first link $\theta_1(t)$ over $10$ time trajectories. The green dashed lines represent the boundary of the safe region.
    }}
	\label{fig_3-link}
\end{figure}

{\color{black}
\paragraph{Further investigations.} We further compare the FESSNC with the existing learning based controllers, including SYNC~\cite{zhang2023sync}, NNDM controller (NNDMC)~\cite{mazouzsafety} and RSM controller (RSMC)~\cite{lechner2022stability}, using the kinetic bicycle model. The results provided in Table~\ref{table2} demonstrate the superiority of the FESSNC over the existing methods in terms of stability, safety, scalability and control energy. For the stability guarantee, we achieve the exponential stability (ES), which is stronger than the asymptotic stability (AS). For the safety guarantee, we assure the controlled trajectories stay in the safe region almost surely (a.s.) instead of with some probability $p<1$. We emphasize that we are the first to achieve the rigorously theoretical guarantee for stability and safety, while all the other methods rely on the numerical approximation to some extent in the guarantee process. We provide more details in Appendix~\ref{experiment1_sub}.

\begin{table*}[htb]
     \centering
     \caption{{\color{black} Comparison of learning-based controllers.} }
\label{table2}
%		\resizebox{\linewidth}{!}{
    \begin{tabular}{ccccccc}
				\toprule
				\toprule
				\multicolumn{1}{c}{\multirow{2}{*}{Index}} & \multicolumn{5}{c}{Method} &    \\
				\cmidrule(lr){2-7}  
				 &FESSNC & SYNC & NNDMC & RSMC  & RSMC+ICNN   \\
				\midrule 
				 Training time (sec)	&\textbf{14}	&125&	209	&435& 987	\\
				  Safety rate(\%) &	\textbf{100}&	90&	100	&10&90\\		  
				 Success rate (\%)&	\textbf{100}	&90	&60	&
                    0&90\\
                    Control energy	&\textbf{1.5}&	1.8&	2.6&	14.5&1.8\\
                    Computational Complexity	&$\mathbf{\mathcal{O}(d)}$ &$\mathcal{O}(k^dd^2)$&$\mathcal{O}(k^dd^3)$&$\mathcal{O}(k^d)$&$\mathcal{O}(k^d)$\\
                    Stability guarantee (Type)&	Yes (ES)	&Yes (AS)&	No	&Yes (AS)&Yes (AS)\\
                    Safety guarantee (Type)&	Yes (a.s.)	&Yes (a.s.)&	Yes (probability)	&No&No\\
                    Type of guarantee	&Theoretical&	Numerical	&Numerical&	Numerical&Numerical\\
				\midrule
				\bottomrule
			\end{tabular}
\end{table*}

\paragraph{Extension to the nonparametric settings.} Although we employ the neural controller as the candidate controller, our framework also works for the nonparametric controller by directly applying the approximated operators in Eqs.~\eqref{st proj},\eqref{sa proj} to it. We apply our framework to the kernel machine based controller and compare it with the FESSNC using a $3$-link pendulum control task, the results shown in Figure~\ref{fig_3-link} demonstrate our approximate projection operations have the ability to improve the performance of the kernel based method without further training. More details are provided in Appendix~\ref{experiment 4}}

\subsection{Synchronization of Coupled Oscillators}
\begin{figure}[htp]
	\centering
	\includegraphics[width=8.2cm]{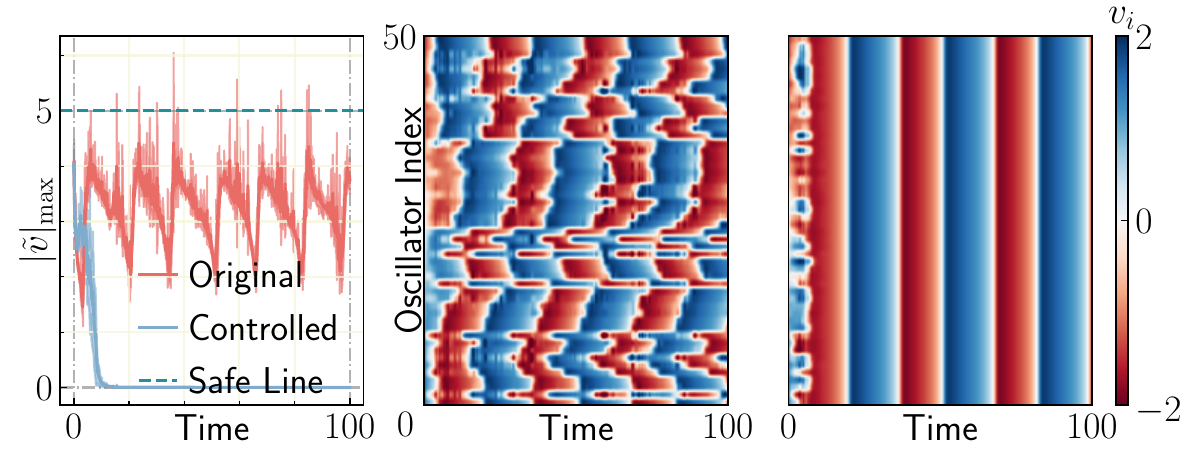}
	\caption{The maximum value of the variances' modulus $|\tilde{v}|_{\max}\triangleq\max_{1\le i\le 50}|v_i|$ along the trajectories of the original dynamics and the controlled dynamics (\textbf{left}). The solid lines are obtained by averaging
		the $5$ sampled trajectories, while the shaded areas stand
		for the variance regions. Dynamics of the small-world network of coupled FHN without control (\textbf{center}) and with control (\textbf{right}).
	 The color represents the membrane potential state $v_i$
		of each oscillator.
	}
	\label{FHN}
 \vskip -0.2in
\end{figure}

The numerical experiments in the previous section demonstrated the feasibility of the proposed framework. Here, we introduce a further
application of FESSNC in the synchronization of coupled oscillators {\color{black}to validate the scalability of FESSNC in high dimensional system}. Phase synchronization phenomenon in coupled oscillatory systems has been extensively studied in the context of collective behavior of complex networks from all scales and domains \citep{rosenblum2001phase,sauseng2008does,strogatz2004sync}. 
In this section, we show that the FESSNC can synchronize the incoherent coupled dynamics in stochastic settings. Specifically, we consider a small world network of $n$ coupled FitzHugh-Nagumo (FHN) models ~\citep{conley1986bifurcation} under the stochastic coupled forces as:
\begin{equation*}
	\begin{aligned}
		&\mathrm{d}v_i=(v_i-\frac{v_i^3}{3}-w_i+1)\mathrm{d}t+\frac{1}{3}\sum_{i=1}^nL_{ij}v_j\mathrm{d}B_t,\\
		&\mathrm{d}w_i=0.1(v_i+0.7-0.8w_i)\mathrm{d}t,
	\end{aligned}
\end{equation*}
where the fast variable $v$ and the slow variable $w$ correspond to the  membrane potential and the recovery variable in an excitable nerve cell, respectively~\citep{keener1998mathematical}. {\color{black}We select $n=50$ s.t. the networked dynamic is a $100$-dimensional system.} We denote the state variables by $\vx=(v,w)$. Although the single FHN model has a stable limit cycle $\Gamma$, the synchronization manifold $\mathcal{M}\triangleq\{\vx_1(t)\subset\Gamma:\vx_i(t)=\vx_1(t),\forall i\}$ of the coupled system is not stable due to the stochasticity. We now stabilize the coupled system to the synchronization manifold $\mathcal{M}$  by stabilizing the equilibrium of the equation of variance $\tilde{\vx}_i=\vx_i-\vx_1$ with the FESSNC. The safe region is $\{\tilde{\vx}_{1:n}:\max_{1\le i\le n}(\tilde{v}_i,\tilde{w}_i)\le5\}$ under the ZBF  $h(\tilde{\vx}_{i:n})=25-\max_{1\le i\le n}(\tilde{v}_i^2,\tilde{w}_i^2)$. The results are shown in Figure~\ref{FHN}. It can be seen that the controlled coupled FHN dynamics converge to the synchronization manifold quickly without crossing the safe region. The experimental details are provided in Appendix.

\section{Related Works}

\paragraph{Stabilization of Nonlinear Systems via Deep Learning}
Recent works have focused on finding the stabilization controllers for dynamical systems through deep learning~\citep{chang2020neural,schlaginhaufen2021learning,dawson2022safe,zhang2022neural}. 
% All these existing methods focus on training the parameterized controllers on the finite dataset according to the stability conditions. Although they may reduce the loss to zero in the training stage, they cannot guarantee that the stability conditions are satisfied on the whole data space as well.
Prior works have also studied how to rigorously satisfy the expected conditions in the deep learning settings, including projecting the learned dynamics to the stable region~\citep{kolter2019learning} and solving a constrained optimization problem at each parameters updating iteration~\citep{chow2019lyapunov}. The drawbacks of these methods are that they either change the equilibria of the original dynamics or cause large computational costs.
\paragraph{Control Barrier Functions for Safety} Recent works have sought to incorporate reciprocal barrier functions, namely the reciprocal of $h(\vx)$, with  suppressing explosion techniques to devise safe controllers~\citep{clark2019control,jankovic2018robust,8796030}. These works impose complicated conditions for the safety, which makes it difficult for us to construct the projection operator. In contrast to the RBFs, ZBFs have \citep{taylor2020learning,taylor2020adaptive} be used to construct the safety-critical control for ODEs with more simple structures. Dawson et al. have utilized the ZBFs in the neural robotic control framework ~\citep{dawson2022safe}. We notice that only the RBFs for SDEs have been intensely studied while the ZBFs for SDEs have not been investigated yet due to the difficulties in theoretical analyses~\citep{clark2019control,clark2021control}.

\section{Conclusion and Future Works}
We have proposed FESSNC, a two-stage framework that learns the neural controllers for stabilizing SDEs with exponential stability and safety guarantees. We present a new theorem on safety criteria for SDEs based on the ZBFs that are more concise than the existing theoretical results. In addition, we develop projection operators to project the learned controllers to the space of exponentially-stable and safe controllers, and we further present approximate projection operators with closed forms to reduce the computational costs. Notably, we prove the existence of both exponentially stable and safe controllers.
We accelerate the training and testing processes with trace estimation methods. The empirical studies show our framework, an offline control policy, performs better the the SOTA online methods GP-MPC and BALSA. In addition to neural controllers, our projection operators can be applied to the nonparametric controllers to improve their performance.

%We believe that our work is a step towards developing learning-based control algorithms that can possess reliable stability and safety, and can have potential applications in stability-critical and safety-critical real-world scenarios. An exciting direction for future work would be to extend the existing model-based offline framework to online learning or model-free cases. In addition, we consider the deterministic control instead of the stochastic control proposed in~\citep{zhang2022neural} because it is difficult to construct an approximate projection operator for stochastic control and we leave it for future work.

Our work can be applied in stability-critical and safety-critical real-world scenarios. An exciting direction for future work would be to extend the existing model-based offline framework to online learning or model-free cases. Moreover, we consider the deterministic control instead of the stochastic control proposed in~\citep{zhang2022neural} because it's difficult to construct an approximate projection operator for stochastic control, and we leave it for future work. Finally, how to extend our current fully actuated control framework to the under actuated control settings appeared in many real-world scenarios is a challenging direction.

\section*{Impact Statements} 
This paper presents work whose goal is to advance the field of Machine Learning. There are many potential societal consequences of our work, none of which we feel must be specifically highlighted here.

% \clearpage
\section*{Acknowledgements}
J. Zhang is supported by the China Scholarship Council (No. 202306100154). Q. Zhu is supported by the China Postdoctoral Science Foundation (No. 2022M720817), by the Shanghai Postdoctoral Excellence Program (No. 2021091), and by the STCSM (Nos. 21511100200, 22ZR1407300, 22dz1200502, and 23YF1402500). W. Lin is supported by the NSFC (Grant No. 11925103), by the STCSM (Grants No. 22JC1402500 and No. 22JC1401402), and by the SMEC (Grant No. 2023ZKZD04). The computational work presented in this article is supported by the CFFF platform of Fudan University.

% % Acknowledgements should only appear in the accepted version.
% \section*{Acknowledgements}

% \textbf{Do not} include acknowledgements in the initial version of
% the paper submitted for blind review.

% In the unusual situation where you want a paper to appear in the
% references without citing it in the main text, use \nocite
%\nocite{langley00}
% \newpage
\bibliography{main}
% \bibliography{sde}
\bibliographystyle{icml2024}

%%%%%%%%%%%%%%%%%%%%%%%%%%%%%%%%%%%%%%%%%%%%%%%%%%%%%%%%%%%%%%%%%%%%%%%%%%%%%%%
%%%%%%%%%%%%%%%%%%%%%%%%%%%%%%%%%%%%%%%%%%%%%%%%%%%%%%%%%%%%%%%%%%%%%%%%%%%%%%%
% DELETE THIS PART. DO NOT PLACE CONTENT AFTER THE REFERENCES!
%%%%%%%%%%%%%%%%%%%%%%%%%%%%%%%%%%%%%%%%%%%%%%%%%%%%%%%%%%%%%%%%%%%%%%%%%%%%%%%
%%%%%%%%%%%%%%%%%%%%%%%%%%%%%%%%%%%%%%%%%%%%%%%%%%%%%%%%%%%%%%%%%%%%%%%%%%%%%%%

\clearpage
\newpage

\appendix

\onecolumn

\section{Appendix}
\subsection{Proofs and Derivations}\label{proofs}
\subsubsection{Notations and Preliminaries}
In this section, we introduce some basic definitions and notations and then provide the proofs of the theoretical results. 
\paragraph{Notations.} Throughout the paper, we employ the following notation. Let $(\Omega,\mathcal{F},\{\mathcal{F}_t\}_{t\ge0},\mathbb{P})$ be a complete probability space with a filtration $\{\mathcal{F}_t\}_{t\ge0}$ satisfying the usual conditions (i.e. it is increasing and right continuous while $\mathcal{F}_0$ contains all $\mathbb{P}$-null sets). If $x,y$ are real numbers, then $x\wedge y$ denotes the minimum of $x$ and $y$, $x\vee y$ denotes the maximum of $x$ and $y$, $x\ll y$ denotes $x$ is much smaller than $b$ and $a\gg b$ vice versa. Let $\langle \vx,\vy\rangle$ be the inner product of vectors $\vx,\vy\in\mathbb{R}^d$. For a second continuous function $f(\vx):\mathbb{R}^d\to\mathbb{R}$, let $\nabla f$ denote the gradient of $f(\vx)$, that is, $\mathcal{H}f$ denote the Hessian matrix of $f$. For the two sets $A,B$, let $A\subset B$ denote that $A$ is covered in $B$. Denote by $\log$ the base $e$ logarithmic function. Denote by $\Vert \cdot \Vert$ the $L^2$-norm for any given vector in $\mathbb{R}^d$.  Denote by $\vert \cdot\vert$ the absolute value of a scalar number or the modulus length of a complex number. For $A=(a_{ij})$, a matrix of dimension $d\times r$, denote by $\|A\|^2_{\rm F}=\sum_{i=1}^{d}\sum_{j=1}^{r}a_{ij}^2$ the Frobenius norm.

\begin{definition}(Martingale)
	The stochastic process $X_t$ on $t\ge0$ is called a martingale (submartingale) on probability space $(\Omega,\mathcal{F},\{\mathcal{F}_t\}_{t\ge0},\mathbb{P})$ with a filtration $\{\mathcal{F}_t\}_{t\ge0}$ satisfying the usual conditions, if the following two conditions hold: $\mathrm{(1)}$ $X_t$ is $\mathcal{F}_t-$measurable for any $t\ge0$; $\mathrm{(2)}$ $\mathbb{E}[X_t|\mathcal{F}_s]=X_s(\mathbb{E}[X_t|\mathcal{F}_s]=X_s)$ for any $t>s\ge0$.
	
\end{definition}

\begin{definition}(Stopping Time)  Given probability space $(\Omega,\mathcal{F},\{\mathcal{F}_t\}_{t\ge0},\mathbb{P})$ and a mapping $\tau:\Omega\to[0,\infty)$, we call $\tau$ an $\{\mathcal{F}_t\}_{t\ge0}$ stopping time, for any $t\ge0$, ${\tau\le t}\in\mathcal{F}_t$,
\end{definition}

\begin{definition}(Local Martingale) The stochastic process $X_t,~t\ge0$ is called a local martingale, if there exists a family of stopping times $\{\tau_n\}_{n\in\mathbb{Z}_{+}}$ such that $\lim_{n\to\infty}\tau_n=\infty,~a.s.$ and $\{X_{t\wedge\tau_n}\}_{n\in\mathbb{Z}_{+}}$ is a martingale.
\end{definition}

\begin{definition}(It$\hat{\textrm{o}}$'s Process)
	Let $B_t$ be a $d$-dimensional Brownian motion on probability space $(\Omega,\mathcal{F},\{\mathcal{F}_t\}_{t\ge0},\mathbb{P})$. A ($d$-dimensional)
	It$\hat{\textrm{o}}$'s process is a stochastic process $X_t$ on $(\Omega,\mathcal{F},\{\mathcal{F}_t\}_{t\ge0},\mathbb{P})$ in the form of
	\begin{equation*}
		X_t=X_0+\int_0^t u(s,w)\mathrm{d}s+\int_0^t\mathrm{d}v(s,w)B_s\left(\Leftrightarrow \mathrm{d}X_t=u(t,w)\mathrm{d}t+v(t,w)\mathrm{d}B_t\right),
	\end{equation*}
	where $u$ and $v$ satisfy the constraints as follows:
	\begin{equation*}
		\begin{aligned}
			&\mathbb{P}\left[\int_0^t\Vert v(s,w)\Vert^2\mathrm{d}s<\infty~ \mathrm{for~all}~t\ge0\right]=1,\\
			&\mathbb{P}\left[\int_0^t\Vert u(s,w)\Vert\mathrm{d}s<\infty~ \mathrm{for~all}~t\ge0\right]=1.
		\end{aligned}
	\end{equation*}
\end{definition}

\begin{definition}(It$\hat{\textrm{o}}$'s Formula)
	Let $X_t$ be a $d$-dimensional It$\hat{\textrm{o}}$'s process given by
	\begin{equation*}
		\mathrm{d}X_t=u\mathrm{d}t+v\mathrm{d}B_t.
	\end{equation*}
	Let $f(t,\vx)\in C^{1,2}([0,\infty)\times\mathbb{R}^d)$. Then, $Y_t=f(t,X_t)$ is an It$\hat{\textrm{o}}$'s process as well, satisfying
	\begin{equation*}
		\mathrm{d}Y_t=\dfrac{\partial h}{\partial t}(t,X_t)\mathrm{d}t+\nabla_{\vx} f(t,X_t)\cdot\mathrm{d}X_t+\dfrac{1}{2}\mathrm{Tr}(\mathrm{d}X_t^\top\mathcal{H}f(t,X_t)\mathrm{d}X_t).
	\end{equation*}
	
\end{definition}

% \begin{lem}(Gronwall's inequality)
% let
% \end{lem}

\subsubsection{Proof of Theorem \ref{thm2}}\label{proof thm2}

To begin with, we check the inequality constraint in $\mathcal{U}_{es}(V,\mathcal{X})$ is satisfied by the projection element, that is 
\begin{equation*}\label{st to check}
	\mathcal{L}_{\vu}V\big|_{\vu=\hat{\pi}_{es}(\vu,\mathcal{U}_{es}(V,\mathcal{X}))}\le cV.
\end{equation*}
From the definition of the derivative operator, we have
\begin{equation*}
	\begin{aligned}
	\mathcal{L}_{\vu}V\big|_{\vu=\hat{\pi}_{es}(\vu,\mathcal{U}_{es}(V,\mathcal{X}))}&=\nabla V\cdot(f+\vu-\dfrac{\max(0,\mathcal{L}_{\vu}V-cV)}{\Vert \nabla V\Vert^2}\cdot\nabla V)+\dfrac{1}{2}\mathrm{Tr}\left[g^\top\mathcal{H}Vg\right]\\
	&=\nabla V\cdot(f+\vu)+\dfrac{1}{2}\mathrm{Tr}\left[g^\top\mathcal{H}Vg\right]-\nabla V\cdot \dfrac{\max(0,\mathcal{L}_{\vu}V-cV)}{\Vert \nabla V\Vert^2}\cdot\nabla V\\
	&=\mathcal{L}_{\vu}V-\max(0,\mathcal{L}_{\vu}V-cV)\le cV.
	\end{aligned}
\end{equation*}
Next, we show the equivalent condition of the Lipschitz continuity of projection element. Notice that $\vu\in\text{Lip}(\mathcal{X})$, then we have 
\begin{equation*}
	\hat{\pi}_{es}(\vu,\mathcal{U}_{es}(V,\mathcal{X}))\in\text{Lip}(\mathcal{X}),\iff\dfrac{\max(0,\mathcal{L}_{\vu}V-cV)}{\Vert \nabla V\Vert^2}\cdot\nabla V\in\text{Lip}(\mathcal{X}).
\end{equation*}
Since $	\Vert\dfrac{\nabla V}{\Vert\nabla V\Vert}\Vert$ is a continuous unit vector, and naturally is Lipschitz continuous, we only need to consider the remaining term $\dfrac{\max(0,\mathcal{L}_{\vu}V-cV)}{\Vert \nabla V\Vert}$. According to the definition, all the functions occured in this term are continuous, so we only need to bound this term to obtain the global Lipschitz continuity, that is 
\begin{equation*}
\dfrac{\max(0,\mathcal{L}_{\vu}V-cV)}{\Vert \nabla V\Vert}\in\text{Lip}(\mathcal{X}),\iff\sup_{\vx\in\mathcal{X}}\dfrac{\max(0,\mathcal{L}_{\vu}V-cV)}{\Vert \nabla V\Vert}<+\infty.
\end{equation*}
When $\mathcal{L}_{\vu}V\le cV$, obviously we have $\max(0,\mathcal{L}_{\vu}V-cV)=0< +\infty$, otherwise, since $V\ge\varepsilon\Vert\vx\Vert^p$ and $c<0$, we have
\begin{equation*}
	\mathcal{L}_{\vu}V-cV\ge\mathcal{L}_{\vu}V-c\varepsilon\Vert\vx\Vert^p\approx\mathcal{O}(\Vert\vx\Vert^p)\to\infty(\Vert\vx\Vert\to\infty).
\end{equation*}
Thus, we have 
\begin{equation*}
	\sup_{\vx\in\mathcal{X}}\dfrac{\max(0,\mathcal{L}_{\vu}V-cV)}{\Vert \nabla V\Vert}<+\infty\iff \sup_{\vx\in\mathcal{X}}\Vert\vx\Vert<+\infty,
\end{equation*}
which completes the proof.

\subsubsection{Proof of Example \ref{ex1}}\label{proof ex1}

	Consider a 1-D SDE as follows:
\begin{equation*}
	\mathrm{d}x=ax\mathrm{d}t+bx\mathrm{d}B_t,
\end{equation*}

Then, according to the Ito formula, we have
\begin{equation*}
	\begin{aligned}
			\mathrm{d}\log x &= \dfrac{\mathrm{d}x}{x}-\dfrac{(\mathrm{d}x)^2}{2x^2}\\
		&=a\mathrm{d}t+b\mathrm{d}B_t-\dfrac{b^2}{2}\mathrm{d}t.
	\end{aligned}
\end{equation*}
Integrate the both sides we have 
\begin{equation*}
	\log x(t)=x(0)+((a-\dfrac{b^2}{2})t+bB(t)).
\end{equation*}
Then we have 
\begin{equation*}
	x(t)=x_0\exp((a-\dfrac{b^2}{2})t+bB(t)),~x_0>0.
\end{equation*}
Then for any $M>0$, we have 
\begin{equation*}
	x(t)>M,\iff x(0)\exp((a-\dfrac{b^2}{2})t+bB(t)),\iff B(t)>\dfrac{1}{b}(\log\dfrac{M}{x_0}-(a-\dfrac{b^2}{2})t).
\end{equation*}
Notice that $B(t)\sim \mathcal{N}(0,t)$ is a Gaussian noise, it naturally holds that
\begin{equation*}
	\mathbb{P}(B(t)>\dfrac{1}{b}(\log\dfrac{M}{x_0}-(a-\dfrac{b^2}{2})t))>0.
\end{equation*}
Next, for the exponential stability, let potential $V=\frac{1}{2}x^2$, then we have 
\begin{equation*}
	\mathcal{L}V=(a+\dfrac{b^2}{2})x^2=2(a+\dfrac{b^2}{2})V.
\end{equation*}
Since  $a<0,~|a|>\frac{1}{2}b^2$ and $x_0>0$, we have $a+\dfrac{b^2}{2}<0$, then according to Theorem~\ref{thm1}, the zero solution is exponential stable. The proof is completed.

\subsubsection{Proof of Theorem \ref{thm3}}\label{proof thm3}
	Notice that each sample path of $\vx(t)$ is continuous and $h(\vx)$ is also continuous. This implies that $h(\vx(t))>0\iff \vx(t)\in\mathrm{int}(\mathcal{C})$. Now we prove $h(\vx(t))>0~a.s.$ with initial $h(\vx(0))>0$, which is equivalent to $\tau=\infty~a.s.$, where stopping time $\tau=\inf\{t\ge0:h(\vx(t))=0\}$. we prove it by contradiction. If $\tau=\infty~a.s.$ was false, then we can find a pair of constants $T>0$ and $M\gg1$ for $\mathbb{P}(B)>0$, where $$B=\{w\in\Omega:\tau<T~\text{and}~\Vert\vx(t)\Vert\le M,\forall 0\le t\le T \}.$$
	But, by the standing hypotheses, there exists a positive constant $K_M$ such that
	$$\alpha(x)\le K_Mx,~\forall \vert x\vert\le\sup_{\Vert\vx\Vert\le M}h(\vx)<\infty.$$
	Then, for $w\in B$ and $t\le T$,
	$$\alpha(h(\vx(t)))\le K_Mh(\vx(t)).$$
	Now, for any $\varepsilon\in(0,h(\vx(0)))$, define the stopping time 
	$$\tau_\varepsilon=\inf\{t\ge0:h(\vx(t))\notin(\varepsilon, h(\vx(0))\}.$$
	By It$\hat{\text{o}}$'s formula,
	\begin{equation*}
		\begin{aligned}
			\mathrm{d}h(\vx(t))&=\mathcal{L}h\mathrm{d}t+\nabla h^\top g\mathrm{d}B_t\\
			&\ge -\alpha(h)\mathrm{d}t+\nabla h^\top G\mathrm{d}B_t.
		\end{aligned}
	\end{equation*}
	Take expectation on both sides with respect to $\tau_\varepsilon$ on set $B$,
	\begin{equation*}
		\begin{aligned}
			\mathbb{E}[h(\vx(\tau_\varepsilon\wedge t))\mathbbm{1}_B]&\ge h(\vx(0))-\int_0^{t}\mathbb{E}[\lambda(h(\vx(\tau_\varepsilon\wedge s)))\mathbbm{1}_B]\mathrm{d}s\\
			&\ge h(\vx(0))-\int_0^{t}\mathbb{E}[K_Mh(\vx(\tau_\varepsilon\wedge s))\mathbbm{1}_B]\mathrm{d}s.
		\end{aligned}
	\end{equation*}
	By Gronwall's inequality, 
	\begin{equation*}
		\begin{aligned}
			\mathbb{E}[h(\vx(\tau_\varepsilon\wedge t))\mathbbm{1}_B]&\ge h(\vx(0))e^{-K_M(\tau_\varepsilon\wedge t)\mathbb{P}(B)}.\\
		\end{aligned}
	\end{equation*}
	Note that for $w\in B$, $\tau_\varepsilon\le T$ and $h(\vx(\tau_\varepsilon))=\varepsilon$. The above inequality, therefore, implies that
	\begin{equation*}
		\mathbb{E}[\varepsilon\mathbb{P}(B)]=\varepsilon\mathbb{P}(B)\ge h(\vx(0))e^{-K_M(\tau_\varepsilon\wedge T)\mathbb{P}(B)}\ge h(\vx(0))e^{-K_M T\mathbb{P}(B)}.
	\end{equation*}
	Letting $\varepsilon\to0$ yields that $0\ge h(\vx(0))e^{-K_M T\mathbb{P}(B)}$, but this contradicts the definition of $B$ and $\mathbb{P}(B)>0$.\\
	The proof is completed.

\subsubsection{Proof of Theorem \ref{thm4}}\label{proof thm4}

Similar to section~\ref{proof thm2}, we only have to check the safety constraint is satisfied by the projection element, that is 
\begin{equation*}
	\mathcal{L}_{\vu}h\big|_{\vu=\hat{\pi}_{sf}(\vu,\mathcal{U}_{sf}(\alpha,\mathcal{C}))}\ge-\alpha(h).
\end{equation*}
From the definitions we have
\begin{equation*}
\begin{aligned}
	\mathcal{L}_{\vu}h\big|_{\vu=\hat{\pi}_{es}(\vu,\mathcal{U}_{es}(V,\mathcal{X}))}&=\nabla h\cdot(f+\vu+\dfrac{\max(0,-\mathcal{L}_{\vu}h-\alpha (h))}{\Vert \nabla h\Vert^2}\cdot\nabla h)+\dfrac{1}{2}\mathrm{Tr}\left[g^\top\mathcal{H}hg\right]\\
	&=\nabla h\cdot(f+\vu)+\dfrac{1}{2}\mathrm{Tr}\left[g^\top\mathcal{H}hg\right]+\nabla h\cdot \dfrac{\max(0,-\mathcal{L}_{\vu}h-\alpha (h))}{\Vert \nabla h\Vert^2}\cdot\nabla h\\
	&=\mathcal{L}_{\vu}h+\max(0,-\mathcal{L}_{\vu}h-\alpha(h))\ge -\alpha(h).
\end{aligned}
\end{equation*}
Notice that the state space now is $\mathcal{C}$, which is bounded, then the proof is completed.

\subsubsection{Proof of Theorem \ref{thm5}}\label{proof thm5}

We prove there exists a controller $\vu$ in $\mathcal{U}_{es}(V,\mathcal{C})\cap\mathcal{U}_{sf}(\alpha,\mathcal{C})$ in two steps.
\paragraph{Step 1.} To begin with, we consider the special case that $\vzero\in\arg\max_{\vx\in\mathcal{C}} h(\vx)$ and $\vzero$ is the only maximum point in $\mathcal{C}$. Notice that 
\begin{equation*}
	\{\vu:\mathcal{L}_{u}h\ge-\alpha(h)+a\}\subset\{\vu:\mathcal{L}_{u}h\ge-\alpha(h)\},~\forall a>0.	
\end{equation*}
Thus, we choose $a=kh(\vzero)$ with $k$ to be defined, then we can get $\vu_0\in\{\vu:\mathcal{L}_{u}h\ge-\alpha(h)+kh(\vzero)\}$, s.t. $\vu_0\in\mathcal{U}_{sf}(\alpha,\mathcal{C})$. We now prove that there exists a potential function $V$ s.t. $\vu_0\in\mathcal{U}_{es}(V,\mathcal{C})$.
Since $\mathcal{C}$ is bounded, the class-$\mathcal{K}$ function $\alpha$ has linear bounds in $\mathcal{C}$ as 
\begin{equation*}
	mx\le\alpha(x)\le qx,~
\end{equation*} 
for some $q>0$. Then for $\vu_0$, the following inequality holds,
\begin{equation*}
	\mathcal{L}_{\vu_0}h\ge-qh+kh(\vzero).
\end{equation*}
Now let $k=q$, we have
\begin{equation*}
	\mathcal{L}_{\vu_0}h\ge q(h(\vzero)-h).
\end{equation*}
Define $V(\vx)=h(\vzero)-h(\vx)\ge0$, due to the continuity on the bounded region $\mathcal{C}$, there exists some positive number $l>0$ s.t. $V(\vx)\ge l\Vert\vx\Vert$ in $\mathcal{C}$. Hence, $V$ is a well-defined  potential function. Notice that 
\begin{equation*}
	\begin{aligned}
			\mathcal{L}_{\vu_0}V&=\nabla V\cdot(f+\vu_0)+\dfrac{1}{2}\mathrm{Tr}\left[g^\top\mathcal{H}hg\right],\\
			\nabla V& =-\nabla h,~~\mathcal{H}V=-\mathcal{H}h.
	\end{aligned}
\end{equation*}
We have 
\begin{equation*}
	\mathcal{L}_{\vu_0}V=-\mathcal{L}_{\vu_0}h\le -q(h(\vzero)-h)=-qV,
\end{equation*}
which means $\vu_0\in\mathcal{U}_{es}(V,\mathcal{C})\cap\mathcal{U}_{sf}(\alpha,\mathcal{C})$.

\paragraph{Step 2.} Now we consider the general case. Since $\vzero\in\mathrm{int}(\mathcal{C})$, and $\{h(\vx)=0\}=\partial \mathcal{C}$, there exists a small neighborhood $O(\vzero,\varepsilon),~\varepsilon>0$ of $\vzero$ s.t.  $O(\vzero,\varepsilon)\subset\mathrm{int}(\mathcal{C})$. Then we can modify $h$ on $O(\vzero,\varepsilon)$ without changing the safe region $\mathcal{C}$. Define the smooth approximation of the Dirac function $\delta_\vzero$ as 
\begin{equation*}
	\lambda(\vx)\triangleq \exp\big(-\dfrac{1}{\Vert\vx\Vert^2-\varepsilon}\big).
\end{equation*}
Since $h$ is bounded on $\mathcal{C}$, we can always choose $M\gg 1$ s.t. 
\begin{equation*}
	\vzero\in\arg\max_{\vx\in\mathcal{C}},\tilde{h},~\tilde{h}=h+M\lambda.
\end{equation*}
According to Step 1, we can choose any $\vu\in\{\mathcal{L}_{\vu}\tilde{h}\ge-\alpha(\tilde{h})+q\tilde{h}(\vzero)\}$, and let $\tilde{V} = \tilde{h}(\vzero)-\tilde{h}$, s.t. $\vu\in\mathcal{U}_{es}(\tilde{V},\mathcal{C})$. Now we only need to show that there exists $\vu\in\{\mathcal{L}_{\vu}\tilde{h}\ge-\alpha(\tilde{h})+q\tilde{h}(\vzero)\}\cap\mathcal{U}_{sf}(\alpha,\mathcal{C})$. Notice that 
\begin{equation*}
	\begin{aligned}
		\mathcal{L}_{\vu}\tilde{h}&=\nabla(h+M\lambda)\cdot(f+u)+\dfrac{1}{2}\mathrm{Tr}\left[g^\top\mathcal{H}(h+M\lambda)g\right]\\
		&=\bigg(\nabla h\cdot(f+u)+\dfrac{1}{2}\mathrm{Tr}\left[g^\top\mathcal{H}hg\right]\bigg)+M\bigg(\nabla\lambda\cdot(f+u)+\dfrac{1}{2}\mathrm{Tr}\left[g^\top\mathcal{H}\lambda g\right]\bigg)\\
		&=\mathcal{L}_{\vu}h+M\mathcal{L}_{\vu}\lambda.
	\end{aligned}
\end{equation*} 
From the median theorem, we have 
\begin{equation*}
	\alpha(\tilde{h})=\alpha(h+M\lambda)=\alpha(h)+\alpha^{'}(\xi)M\lambda\le\alpha(h)+LM\lambda,~\xi\in(0,L),~L=\max_{\vx\in\mathcal{C}}h(\vx).
\end{equation*}
Then we have
\begin{equation*}
	\begin{aligned}
			&\mathcal{L}_{\vu}\tilde{h}\ge-\alpha(\tilde{h})+q\tilde{h}(\vzero)\ge-\alpha(h)-LM\lambda+q\tilde{h}(\vzero),\\
		&\mathcal{L}_{\vu}\tilde{h}\ge-\alpha(h)-LM\lambda+q\tilde{h}(\vzero),\iff\mathcal{L}_{\vu}h\ge-\alpha(h)-LM\lambda+q\tilde{h}(\vzero)-M\mathcal{L}_{\vu}\lambda.
	\end{aligned}
\end{equation*}
Let $N=\max(0,\sup_{\vx\in O(\vzero,\varepsilon)}-LM\lambda+q\tilde{h}(\vzero)-M\mathcal{L}_{\vu}\lambda)$, then we have 
\begin{equation*}
	\mathcal{L}_{\vu}h\ge-\alpha(h)-LM\lambda+q\tilde{h}(\vzero)-M\mathcal{L}_{\vu}\lambda\ge-\alpha(h)+N.
\end{equation*}
From step 1, we know that $\{\mathcal{L}_{\vu}h\ge-\alpha(h)+N\}\subset\mathcal{U}_{sf}(\alpha,\mathcal{C})$, now we choose $\vu\in\{\mathcal{L}_{\vu}h\ge-\alpha(h)+N\}$ and complete the proof.

{\color{black}
For the relaxation in Remark~\ref{rem1}, we notice that when safe region is covered in unbounded cube as that in Remark~\ref{rem1}, the safe function $h$ only depends on the variables that locate in the bounded interval. Without loss of generality, we set the first $j$ variables locate in bounded safe region. Then $h$ only depends on $\vx=(x_1,\cdots,x_j)^\top$ and is bounded in $[a_1,b_1]\times\cdots\times[a_j,b_j]$, hence there still exist $m,~q$ s.t. $mx\alpha(x)\le qx$ as required in \textbf{Step 1}. Then we repeat the rest of the above proof to complete the proof of Remark~\ref{rem1}.
}

\section{Experimental Configurations}\label{appen_details}
In this section, we provide the detailed descriptions for the experimental configurations of the control problems in the main text.  The computing device that we use for calculating our examples includes a single i7-10870 CPU with 16GB memory, and we train all the parameters with Adam optimizer.   
{Our code is available at  \href{https://github.com/jingddong-zhang/FESSNC}{\texttt{https://github.com/jingddong-zhang/FESSNC}}.}

\subsection{Neural Network Structures}

\begin{itemize}
	\item For constructing the potential function $V$, we utilize the ICNN as~\citep{icnn}:
	\begin{equation*}
		\begin{aligned}
			{\vz}_1 &= \sigma(W_0{\vx}+b_0),\\
			{\vz}_{i+1} &= \sigma(U_i{\vz}_i+W_i{\vx}+b_i),\ i=1,\cdots,k-1,\\
			p({\vx}) &\equiv {\vz}_k,\\
			V({\vx}) &= \sigma(p({\vx})-p(\vzero))+\varepsilon\Vert {\vx}\Vert^2,\\ 
		\end{aligned}
	\end{equation*}
	where $\sigma$ is the smoothed \textbf{ReLU} function as defined in the main text, $W_i\in\mathbb{R}^{h_i\times d},\ U_i\in(\mathbb{R}_{+}\cup \{0\})^{h_i\times h_{i-1}}$, $\vx\in\mathbb{R}^d$, and, for simplicity, this ICNN function is denoted by ICNN$(h_0,h_1,\cdots,h_{k-1})$; 
	
	\item The class-$\mathcal{K}$ function $\alpha$ is constructed as:
	\begin{equation*}
		\begin{aligned}
			\vq_1&=\text{ReLU}(W_0s+b_1),\\
			\vq_{i+1}&=\text{ReLU}(W_i\vq_i+b_i),i=1,\cdots,k-2,\\
			\vq_k&=\text{ELU}(W_{k-1}\vq_{k-1}+b_{k-1}),\\
			\alpha(x)&=\int_0^xq_k(s)\mathrm{d}s
		\end{aligned}
	\end{equation*}
	where $W_i\in\mathbb{R}^{h_{i+1}\times h_i}$, and this class-$\mathcal{K}$ function is denoted by $\mathcal{K}(h_0,h_1,\cdots,h_{k})$;
	
	\item The neural control function (nonlinear version) is constructed as: \\ 
	\begin{equation*}
		\begin{aligned}
			\vz_1 &= \mathcal{F}(W_0\vx+B_1),\\
			\vz_{i+1} &= \mathcal{F}(W_i\vz_{i}+b_i),\ i=1,\cdots,k-1,\\
			\textbf{NN}(\vx) &\equiv W_k\vz_{k},\\
			\vu(\vx) &= \text{diag}(\vx)\textbf{NN}(\vx),
		\end{aligned}
	\end{equation*}
	where $\mathcal{F}(\cdot)$ is the activation function,  $W_i\in\mathbb{R}^{h_{i+1}\times h_i}$, and this control function is denoted by Control$(h_0,h_1,\cdots,h_{k+1})$.
\end{itemize}

\subsection{Double Pendulum}\label{experiment1}
Here we model the state of the fully actuated noise-perturbed double pendulum as 
\begin{equation*}
	\begin{aligned}
		&\mathrm{d}\theta_1=z_1\mathrm{d}t\\
		&\mathrm{d}z_1=\dfrac{m_2g\sin(\theta_2)\cos(\theta_1-\theta_2)-m_2\sin(\theta_1-\theta_2)[l_1z_1^2\cos(\theta_1-\theta_2)+l_2z_2^2]-(m_1+m_2)g\sin(\theta_1)}{l_1[m_1+m_2\sin^2(\theta_1-\theta_2)]}\mathrm{d}t+\sin(\theta_1)\mathrm{d}B_t\\
		&\mathrm{d}\theta_2(t)=z_2\mathrm{d}t,\\
		&\mathrm{d}z_2=\dfrac{(m_1+m_2)[l_1z_1^2\sin(\theta_1-\theta_2)-g\sin(\theta_2)+g\sin(\theta_1)\cos(\theta_1-\theta_2)]+m_2l_2z_2^2\sin(\theta_1-\theta_2)\cos(\theta_1-\theta_2)}{l_2[m_1+m_2\sin^2(\theta_1-\theta_2)]}\mathrm{d}t+\sin(\theta_2)\mathrm{d}B_t.
	\end{aligned}
\end{equation*}
We transform the equilibrium to the zero by the coordinate transformation $\tilde{\theta}_{1,2}=\theta_{1,2}-\pi$. the The ZBF is $h(\theta_{1,2},z_{1,2})=\sin(\theta_1)+0.5\iff h(\tilde{\theta}_{1,2},z_{1,2})=0.5-\sin(\tilde{\theta}_1)$. Then we consider the dynamic for $\vx=(\tilde{\theta}_{1,2},z_{1,2})$ For training controller $\vu$, we sample $500$ data from the safe region $[-\pi/6-\pi,7\pi/6-\pi]\times[-5,5]^3$, we construct the NNs as follows.

\paragraph{FESSNC} We parameterize $V(\vx)$ as ICNN$(4,12,12,1)$, $\alpha(x)$ as $\mathcal{K}(1,10,10,1)$, $\vu(\vx)$ as Control$(4,12,12,4)$. We set $\varepsilon=1\mathrm{e}$-$3$, $c=-0.1$. We train the parameters with $\mathrm{lr}=0.1$ for $300$ steps.

\paragraph{GP-MPC} Since we have known the dynamics, there is no need to fitting the dynamics as GP flows from data. We obtain the dynamics of moment $\vz=(\mu,\Sigma)$ of distribution $p(\vx)$ as 
\begin{equation*}
	\begin{aligned}
		\mu_{t+1}=&\mu_t+dt(\cdot f(\mu_t)+\vu_t)\\
		\Sigma_{t+1}=&\Sigma_t+dt\cdot g(\mu_t)
	\end{aligned}
\end{equation*}
where $f,g$ are the drift and diffusion terms, respectively. We denote the above dynamic as 
\begin{equation*}
	\vz_{t+1}=\tilde{f}(\vz_t,\vu_t).
\end{equation*}
Then we solve the MPC problem as 
\begin{equation*}
	\begin{aligned}
		\min_{\vu_{0:k-1}}&~\Phi(\vz_k)+\sum_{i=0}^{k-1}l(\vz_i,u_i),\\
			&\text{s.t.}~\vz_{i+1}=\tilde{f}(\vz_i,\vu_i),~\vz_0=(\vx,\vzero),
		\end{aligned}
	\end{equation*}
	and apply the first control element $\vu_0$ to the current dynamic. We use the PMP method as proposed in~\citep{kamthe2018data} to solve the above problem. For safety and stability, we select cost function as 
	\begin{equation*}
		\begin{aligned}
			\Phi(\vz)=&p_1\Vert\mu\Vert^2+p_3\Vert\Sigma\Vert^2,\\
			l(\vz,\vu)=&p_2\text{ReLU}(\sin(\tilde{\theta}_1)-0.5)+\Vert\vz\Vert^2.
		\end{aligned}
	\end{equation*}
	We set $p_1=6,p_2=2.5,p_3=0.5,k=10$ and find the optimal control sequence via the gradient descent of the Hamiltonian function in PMP.
	
	\paragraph{QP} The object function in QP is set as 
	
	\begin{equation*}
		\begin{aligned}
			\min_{\vu,d_1,d_2}&\frac{1}{2}\Vert\vu\Vert^2+p_1d_1^2+p_2d_2^2,\\
			\text{s.t.}&\mathcal{L}_\vu V-\varepsilon V\le d_1,\\
			&\mathcal{L}_{\vu}\frac{1}{h(\vx)}-\gamma h(\vx)\le d_2,
		\end{aligned}
	\end{equation*}
	where $d_{1,2}$ are the relaxation numbers. We choose $V=\frac{1}{2}\Vert\vx\Vert^2,p_1=10,p_2=10,\varepsilon=0.5,\gamma=2$.

	We use Euler–Maruyama numerical scheme to simulate the system without and with control, and the random seeds are set as $\{1,4,6,8,9\}$.

\subsection{Planar kinematic bicycle model}\label{experiment2}
Here we model the state of the common noise-perturbed kinematic bicycle system as $\vx=(x,y,\theta,v)^\top$, where $x,y$ are the coordinate positions in phase plane, $\theta$ is the heading, $v$ is the velocity. The dynamic is as follows,
\begin{equation*}
	\begin{aligned}
		&\mathrm{d}x(t)=v(t)\cos\theta(t)\mathrm{d}t+x(t)\mathrm{d}B_t,\\
		&\mathrm{d}y(t)=v(t)\sin\theta(t)\mathrm{d}t+y(t)\mathrm{d}B_t,\\
		&\mathrm{d}\theta(t)=v(t)\mathrm{d}t,\\
		&\mathrm{d}v(t)=(x(t)^2+y(t)^2)\mathrm{d}t.
	\end{aligned}
\end{equation*}
The ZBF is $h(x,y,\theta,v)=2^2-(x^2+y^2)$. For training controller $\vu$, we sample $500$ data $(r\cos(w),r\sin(w),\theta,v)$ from the safe region $(r,w,\theta,v)\in[0,3]\times[0,2\pi]\times[-3,3]^2$, we construct the NNs as follows.

\paragraph{FESSNC} We parameterize $V(\vx)$ as ICNN$(4,12,12,1)$, $\alpha(x)$ as $\mathcal{K}(1,10,10,1)$, $\vu(\vx)$ as Control$(4,12,12,4)$. We set $\varepsilon=1\mathrm{e}$-$3$, $c=-0.5$. We train the parameters with $\mathrm{lr}=0.05$ for $500$ steps.

\paragraph{GP-MPC} Since we have known the dynamics, there is no need to fitting the dynamics as GP flows from data. We obtain the dynamics of moment $\vz=(\mu,\Sigma)$ of distribution $p(\vx)$ as 
\begin{equation*}
	\begin{aligned}
		\mu_{t+1}=&\mu_t+dt(\cdot f(\mu_t)+\vu_t),\\
		\Sigma_{t+1}=&\Sigma_t+dt\cdot g(\mu_t),
	\end{aligned}
\end{equation*}
 where $f,g$ are the drift and diffusion terms, respectively. We denote the above dynamic as 
 \begin{equation*}
 	\vz_{t+1}=\tilde{f}(\vz_t,\vu_t)
 \end{equation*}
 Than we solve the MPC problem as 
 \begin{equation*}
 	\begin{aligned}
 		 	\min_{\vu_{0:k-1}}&~\Phi(\vz_k)+\sum_{i=0}^{k-1}l(\vz_i,u_i),\\
 		 	&\text{s.t.}~\vz_{i+1}=\tilde{f}(\vz_i,\vu_i),~\vz_0=(\vx,\vzero),
 	\end{aligned}
 \end{equation*}
and apply the first control element $\vu_0$ to the current dynamic. We use the PMP method as proposed in~\citep{kamthe2018data} to solve the above problem. For safety and stability, we select cost function as 
\begin{equation*}
	\begin{aligned}
		\Phi(\vz)=&p_1\Vert\mu\Vert^2+p_3\Vert\Sigma\Vert^2,\\
		l(\vz,\vu)=&p_2\text{ReLU}(\vz_x^2+\vz_y^2-2^2)+\Vert\vz\Vert^2.
	\end{aligned}
\end{equation*}
We set $p_1=5,p_2=2.5,p_3=0.5,k=5$ and find the optimal control sequence via the gradient descent of the Hamiltonian function in PMP.

\paragraph{QP} The object function in QP is set as 

\begin{equation*}
	\begin{aligned}
	\min_{\vu,d_1,d_2}&\frac{1}{2}\Vert\vu\Vert^2+p_1d_1^2+p_2d_2^2,\\
	\text{s.t.}&\mathcal{L}_\vu V-\varepsilon V\le d_1,\\
	&\mathcal{L}_{\vu}\frac{1}{h(\vx)}-\gamma h(\vx)\le d_2,
	\end{aligned}
\end{equation*}
where $d_{1,2}$ are the relaxation numbers. We choose $V=\frac{1}{2}\Vert\vx\Vert^2,p_1=10,p_2=10,\varepsilon=0.5,\gamma=2$.

We use Euler–Maruyama numerical scheme to simulate the system without and with control, and the random seeds are set as $\{3,5,6,9,10\}$.

{\color{black}
\subsubsection{Further comparison with learning-based controllers}\label{experiment1_sub}
We further compare the proposed FESSNC with the neural network based control methods to stabilize SDEs, we consider the SYNC~\cite{zhang2023sync}, NNDMC~\cite{mazouzsafety} and RSMC~\cite{lechner2022stability} methods. All these considered methods provide stability or safety guarantee for the controller based on the finite decomposition of the compact state space and a fine design of the loss function, which belong to the numerical guarantee and cannot assure the obtained controllers satisfy the stability or safety guarantee rigorously. The model details of the compared methods are listed as follows,
\paragraph{SYNC} We use the grid discretization method with mesh size $r=10$ to obtain $10^4$ grid data in safe region $[-2,2]^2\times[-3,3]^2$. We parameterize the controller as Control$(4,12,12,4)$, the batch size is set as $N=500$, we train the parameters with $\mathrm{lr}=0.05$ for $500$ steps with the same loss function as follows
\begin{equation*}
\begin{aligned}
    L_{sf}({\vtheta},\vtheta_\alpha) &= \dfrac{1}{N}\sum_{i=1}^N\max\left\{0,-\mathcal{L}h(\vx)-\alpha(h(\vx))+4Mr\right\}.\\
 L_{st}({\vtheta}) &=
			\frac{1}{N}\sum_{i=1}^N
			\big[\max\big(0,(\alpha-2)\|{\vx_i}^{\top}g({\vx_i})\|^2
			+\|{\vx}_i\|^2(2\langle {\vx}_i,f_\vu({\vx}_i)\rangle+
			\|g({\vx}_i)\|_{\rm F}^2)\big)\big],\\
   L &= L_{sf}+L_{st}
\end{aligned}
\end{equation*}

Here $f_\vu$ and $g$ represent the controlled drift and diffusion terms, respectively. We set the stability strength as $\alpha=0.8$, the Lipschitz constant of the corresponding functions in the safe region can be calculated as $M=4\max q(x)$ where $\alpha(x)=\int_0^xq(z)\mathrm{d}z$.

\paragraph{NNDMC} We first train the stabilization controller $\vu_{st}$ with the same structure as that in FESSNC, then, we obtain the neural network dynamical model (NNDM) under control. In order to provide safety guarantee for the NNDM, we discretize the safe region to grid data as same in SYNC, then we use the 
\href{https://github.com/Verified-Intelligence/auto_LiRPA}{\texttt{auto LiRPA}} package to compute the linear lower bound $lb_q=(lb_{qx},lb_{qy},lb_{q\theta},lb_{qv})^\top$ and upper bound $ub_q=(ub_{qx},ub_{qy},ub_{q\theta},ub_{qv})^\top$ of the NNDM on each grid subspace $q\in Q$. Here $Q$ is the set of all grid subspace. The input to the auto LiRPA model is the center point of the grid subspace, i.e., the center point of $(x,y,\theta,v)\in[a_1,b_1]\times[a_2,b_2]\times[a_3,b_3]\times[a_4,b_4]$ is $(\frac{a_1+b_1}{2},\frac{a_2+b_2}{2},\frac{a_3+b_3}{2},\frac{a_4+b_4}{2})$. Since we know the safe function is this case, we do not have to solve the SOS problem to obtain the safe function and the minimizer of the safe function is $(0,0,0,0)$. We then use the obtained lower and upper bound to solve the linear programming (LP) problem to obtain the safe controller $\vu_{sf}$. Finally, the controller $\vu=\vu_{sf}+\vu_{st}$ is applied to the original system. According to \cite{mazouzsafety}, the LP problem is solved as follows,
\begin{equation*}
    \begin{aligned}
        \min_{a_1,a_2,a_3,a_4,u_{qx},u_{qy},u_{q\theta},u_{qv}}&\sum_{i=1}^4a_i\\
        \text{s.t.}~&-a_1\le x\le a_1,\\
        &-a_2\le y\le a_2,\\
        &-a_3\le \theta\le a_3,\\
        &-a_4\le v\le a_4,\\
        &lb_{qx}+u_{qx}\le x\le ub_{qx}+u_{qx},\\
        &lb_{qy}+u_{qy}\le y\le ub_{qy}+u_{qy},\\
        &lb_{q\theta}+u_{q\theta}\le \theta\le ub_{q\theta}+u_{q\theta},\\
        &lb_{qv}+u_{qv}\le v\le ub_{qv}+u_{qv},~\forall q\in Q.\\
    \end{aligned}
\end{equation*}
The final safe controller takes value as $\vu(x,y,\theta,v)=(u_{qx},u_{qy},u_{q\theta},u_{qv}),~\forall (x,y,\theta,v)^\top\in q$.

\paragraph{RSMC} According to \cite{lechner2022stability}, the training data should be discretization of the compact state space, here we employ the same grid data as above. For the supermartingale function $V$, we parameterize it with the standard forward neural network FNN$(4,12,12,1)$. The training parameters are the same as above. The loss function is set as 
\begin{equation*}
\begin{aligned}
       L_{st} &= \dfrac{1}{N}\sum_{i=1}^j\left[\text{ReLU}\left(\dfrac{1}{m}\sum_{k=1}^mV(\vz_i^{k})-V(\vz_i)+\tau K\right)\right]\\
       \vz_i^k&=\vz_i+dt\cdot f_\vu(\vz_i)+\sqrt{dt}\xi_kg(\vz_i),~\vz=(x,y,\theta,v)^\top,~\xi_k\sim\mathcal{N}(0,I),\\
       L_{{lip}}& = \text{ReLU}\left(\text{Lip}_V-\dfrac{\kappa}{\tau(\text{Lip}_f(\text{Lip}_\vu+1)+1)}\right),\\
       L &=L_{st}+L_{lip}.
       \end{aligned}
\end{equation*}
We set the tolerance error $\kappa=0.1$, $\tau=6/10$, $m=5$ is the sample size to estimate the expectation, the Euler step size $dt=0.01$, and the Lipschitz constants $\text{Lip}_V$, $\text{Lip}_\vu$ and calculated by the method in \cite{goodfellow2014explaining}.

\paragraph{RSMC+ICNN} We directly replace the FNN supermartingale $V$ as ICNN$(4,12,12,1)$.

We use Euler–Maruyama numerical scheme to simulate the system without and with control, and the random seeds are set as $\{3,6,9,10,11,12,14,15,16,28\}$. The results are summarized in the Table~\ref{table1}. We train all the methods with $\texttt{batch~size}=500$ for $500$ steps and calculate the training time. We test the learned controller on $10$ sample temporal trajectories over $20$ seconds. The results are summarized in the following table.  The safety rate is the ratio of time that the controlled state stays in the safe region to total time. The success rate is the number of controlled states that satisfy the safety constraint and is closer than 0.1 to the target position for consecutive 2 seconds to $10$. The control energy is calculated as $\int_0^T||u(t)||^2dt$.  We notice that the RSMC failed in the task, but we could improve the performance of this method if we replace the original parameterized $V$ function in [R3] with the ICNN constructed in our paper. We can see that our FESSNC achieves the best performance with rather low energy. For the computational complexity, $d$ is the state dimension and $k$ is the grid size for discretizing the state space needed in all the methods in [R1-R3]. The results show that only our FESSNC avoids the curse of dimensionality, and hence our FESSNC can scale to high-dimensional tasks. For the stability guarantee, only our method achieves exponential stability (ES) while the other methods achieve weaker asymptotic stability (AS). For the safety guarantee, our method holds almost surely (a.s.) for any controlled trajectory, that is $\mathbb{P}(x(t)\in \mathcal{C}, t\ge0)=1$, while the NNDMC method holds in probability, that is $\mathbb{P}(x(t)\in \mathcal{C}, 0\le t\le N)>1-N\delta$ for some $\delta\in(0,1)$. Last but not least, all the existing methods provide the numerical guarantee for the neural controller using the finite decomposition of the compact state space as the training data, but none of these methods can assure the stability or safety condition is rigorously satisfied after training. In contrast, we are the first to provide the theoretical stability and safety guarantee for the learned controller based on the analytical approximate projections.

\begin{table}[htb]
     \centering
     \caption{{\color{black} Performance and comparison of learning-based controllers.} }
\label{table1}
%		\resizebox{\linewidth}{!}{
    \begin{tabular}{ccccccc}
				\toprule
				\toprule
				\multicolumn{1}{c}{\multirow{2}{*}{Index}} & \multicolumn{5}{c}{Method} &    \\
				\cmidrule(lr){2-7}  
				 &FESSNC & SYNC & NNDMC & RSMC  & RSMC+ICNN   \\
				\midrule 
				 Training time (sec)	&\textbf{14}	&125&	209	&435& 987	\\
				  Safety rate(\%) &	\textbf{100}&	90&	100	&10&90\\		  
				 Success rate (\%)&	\textbf{100}	&90	&60	&
                    0&90\\
                    Control energy	&\textbf{1.5}&	1.8&	2.6&	14.5&1.8\\
                    Complexity	&$\mathbf{\mathcal{O}(d)}$ &$\mathcal{O}(k^dd^2)$&$\mathcal{O}(k^dd^3)$&$\mathcal{O}(k^d)$&$\mathcal{O}(k^d)$\\
                    Stability guarantee (Type)&	Yes (ES)	&Yes (AS)&	No	&Yes (AS)&Yes (AS)\\
                    Safety guarantee (Type)&	Yes (a.s.)	&Yes (a.s.)&	Yes (probability)	&No&No\\
                    Type of guarantee	&Theoretical&	Numerical	&Numerical&	Numerical&Numerical\\
				\midrule
				\bottomrule
			\end{tabular}
\end{table}
  
\subsubsection{Ablation study}
We further investigate the influence of the hyperparameters in our framework, including $k,\lambda_1,\lambda_2$ relating to the strength of exponential stability, the weight of stability loss, and the weight of safety loss, respectively. We train the FESSNC under different combinations of hyperparameters $\{k,\lambda_1,\lambda_2\}\in\{0.1,0.5,1.0\}^3$, then we test them using $10$ sample trajectories. We summarize the results in the following table, in which the Error is the minimum distance between the controlled bicycle to the target position, i.e., $\text{Error}=\min_t r(t)=||(x(t),y(t))||$, the Std is the standard deviation of $r(t)$, and we use $\max r(t)$ to assess the risk to cross the boundary of safe region $r\le2$. The controllers under all the hyperparameter combinations perform fairly well, implying our method is not very sensitive to the selection of the hyperparameters. 

\begin{table}[htb]
     \centering
     \caption{{\color{black} Performance and comparison of learning-based controllers.} }
\label{table1}
%		\resizebox{\linewidth}{!}{
    \begin{tabular}{ccccccccccc}
				\toprule
				\toprule
				\multicolumn{1}{c}{\multirow{2}{*}{}} & \multicolumn{1}{c}{\multirow{2}{*}{}} & \multicolumn{3}{c}{$\lambda_1=0.1$} &  \multicolumn{3}{c}{$\lambda_1=0.5$}&\multicolumn{3}{c}{$\lambda_1=1.0$}  \\
				\cmidrule(lr){3-11}  
				& &$k=0.1$ & $k=0.5$ & $k=1.0$ & $k=0.1$  & $k=0.5$ & $k=1.0$ & $k=0.1$ & $k=0.5$ & $k=1.0$   \\
				\midrule 
				 \multicolumn{1}{c}{\multirow{3}{*}{Error}} &$\lambda_2=0.1$&$0.02$&$0.02$&$0.02$&\textbf{4e-4}&$0.01$&$0.05$&$0.018$&$0.05$&$0.01$\\		  
	&$\lambda_2=0.5$&$0.11$&$0.03$&$0.03$&\textbf{2e-3}&$0.02$&$0.01$&$0.03$&$0.05$&$0.04$\\
 &$\lambda_2=1.0$&$0.08$&$0.04$&$0.03$&\textbf{8e-4}&$0.01$&$0.06$&$0.03$&$0.01$&$0.07$\\
 \cdashline{1-11}[1pt/1pt]
  \multicolumn{1}{c}{\multirow{3}{*}{Std}} &$\lambda_2=0.1$&$0.06$&$0.38$&$0.35$&\textbf{1e-3}&$0.41$&$2.07$&$0.20$&$0.27$&$0.10$\\		  
&$\lambda_2=0.5$&$0.29$&$0.32$&$0.34$&\textbf{2e-3}&$0.04$&$0.05$&$0.17$&$0.68$&$0.10$\\
 &$\lambda_2=1.0$&$0.29$&$0.58$&$0.38$&\textbf{4e-3}&$0.26$&$0.16$&$0.33$&$0.15$&$0.21$\\
  \cdashline{1-11}[1pt/1pt]
    \multicolumn{1}{c}{\multirow{3}{*}{$\max r(t)$}} &$\lambda_2=0.1$&$1.81$&\textbf{1.73}&$1.90$&$1.94$&$1.82$&$1.83$&$1.78$&$1.95$&$1.83$\\		  
&$\lambda_2=0.5$&$1.91$&\textbf{1.66}&$1.86$&$1.90$&|$1.77$&$1.83$&$1.76$&$1.87$&$1.81$\\
&$\lambda_2=1.0$&$1.72$&$1.70$&$1.78$&$1.88$&$1.89$&$1.77$&\textbf{1.67}&$1.80$& $1.77$\\
				\midrule
				\bottomrule
			\end{tabular}
\end{table}

}

\subsection{Fitzhugh-Nagumo model}\label{experiment3}
Consider the coupled FHN model $\vx_{1:50}$ as 
\begin{equation*}
	\mathrm{d}\vx_i=f(\vx_i)\mathrm{d}t+\sum_{j=1}^Ng(\vx_j)\mathrm{d}B_t
\end{equation*}
where 
\begin{equation*}
	\begin{aligned}
		\vx_i&=v_i,w_i,\\
		f(\vx_i)&=(v_i-\frac{v_i^3}{3}-w_i+1,0.1(v_i+0.7-0.8w_i))^\top,\\
		g(\vx_i)&=(\frac{1}{3}v_i,0)^\top
	\end{aligned}
\end{equation*}
Let $\tilde{\vx}=(\tilde{\vx}_1,\cdots,\tilde{\vx}_{50})^\top=(\vx_1-\vs,\cdots,\vx_{50}-\vs)^\top$, where $\vs$ is on the synchronization manifold s.t. $\dot{\vs}=f(\vs)$, then we have the following variance equation
\begin{equation*}
	\mathrm{d}\delta=(I\otimes \nabla f)\delta\mathrm{d}t+(I\otimes\nabla g)\delta\mathrm{d}B_t.
\end{equation*}
We only need to stabilize the variance $\delta$ to zero, and the ZBF is set as $h(\delta)=5^2-\max_{1\le i\le50}(\tilde{v}_i^2,\tilde{w}_i^2)$. We parameterize $V(\vx)$ as ICNN$(100,100,100,1)$, $\alpha(x)$ as $\mathcal{K}(1,10,10,1)$, $\vu(\vx)$ as Control$(100,200,200,100)$. We set $\varepsilon=1\mathrm{e}$-$3$, $c=-0.1$. We train the parameters with $\mathrm{lr}=0.1$ for $300$ steps. We use Euler–Maruyama numerical scheme to simulate the system without and with control, and the random seeds are set as $\{1,4,5,9,15\}$.

{\color{black}
\subsection{3-link Pendulum}\label{experiment 4}
In this section, we test our framework on a complex $3$-link planar pendulum, which possesses $6$ state variables $(\theta_1,\theta_2,\theta_3,\dot{\theta}_1,\dot{\theta}_2,\dot{\theta}_3)$ representing the $3$ link angles and the $3$ angle velocities of three joints. For sake of simplicity, we suppose the joints and links have unit mass $m_i$, unit length $l_i$, unit moment of inertia $I_i$, and unit relative position $l_{ci}$ of the center of gravity. The corresponding controlled system under noise perturbed is:
	\begin{equation*}
	    \begin{aligned}
	     &\mathrm{d}\vx=\vy\mathrm{d}t,\\
	    &\mathrm{d}\vy=\left\{\bm M(\vx)^{-1}[-\bm N(\vx,\vy)\vy-\bm Q(\vx)]+\vu(\vx,\vy)\right\}\mathrm{d}t+\vg(\vx)\mathrm{d}B_2(t),\\
	    \end{aligned}
	\end{equation*}
where $\vx=(\theta_1,\theta_2,\theta_3)^\top$, $\vy=(\dot{\theta}_1,\dot{\theta}_2,\dot{\theta}_3)^\top$, $\bm M,\bm N\in\mathbb{R}^{3\times3}$, $\bm Q\in\mathbb{R}^3$, $\vg(\vx)=(\sin(\theta_1),\sin(\theta_2),\sin(\theta_3))^\top$ with 
$$
\begin{array}{l}
  \displaystyle 
  M_{ij}=a_{ij}\cos(x_j-x_i),~
  ~ N_{ij}=-a_{ij}y_j\sin(x_j-x_i),~
  ~ Q_i=-b_i\sin(x_i),
  \\
 a_{ii}=I_i+m_il_{ci}^2+l_i^2\sum_{k=i+1}^{3}m_k, 
  ~~a_{ij}=a_{ji}=m_jl_il_{cj}+l_il_j\sum_{k=j+1}^3m_k,\\
 b_i=m_il_{ci}+l_i\sum_{k=i+1}^3m_k.
 \end{array}
$$ 
\begin{figure}[htp]
	\centering
	\includegraphics[width=15cm]{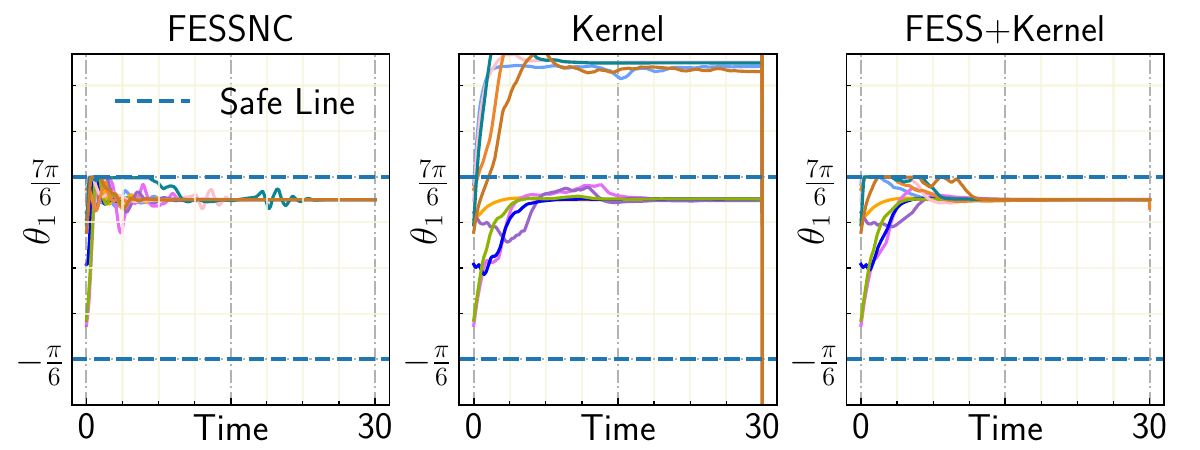}
	\caption{{\color{black}The angle of the first link $\theta_1(t)$ over $10$ time trajectories. The green dashed lines represent the boundary of the safe region.
    }}
	\label{3-link}
\end{figure} 

The mission is to stabilize the 3-link pendulum to the upright position $\theta_i=\pi,~i=1,2,3$ without crossing the safe region $\theta_1\in(-\frac{\pi}{6},\frac{7\pi}{6})$ (the same safe region as that for double pendulum in Section~\ref{experiment1}). In addition to test the efficacy the proposed FESSNC, we also investigate the extension of our framework to nonparametric setting. We employ the kernel machine based controller as a benchmark of nonparametric controller. Specifically, we consider the implementation in Rectified Flow~\cite{liu2022flow} with Gaussian RBF kernel to design a model-based stabilization controller, named Kernel. This controller naturally lacks stability and safety guarantee, hence we apply our approximate projection operators to this controller to provide expected guarantees, the corresponding controller is dubbed as FESS+Kernel.  The detailed model structures are set as follows,

We transform the equilibrium to the zero by the coordinate transformation $\tilde{\theta}_{1,2}=\theta_{1,2}-\pi$. the The ZBF is $h(\theta_{1,2},z_{1,2})=\sin(\theta_1)+0.5\iff h(\tilde{\theta}_{1,2},z_{1,2})=0.5-\sin(\tilde{\theta}_1)$. Then we consider the dynamic for $\vx=(\tilde{\theta}_{1,2},z_{1,2})$. For training controller $\vu$, we sample $500$ data from the safe region $[-\pi/6-\pi,7\pi/6-\pi]\times[-5,5]^3$, we construct the NNs as follows.

\paragraph{FESSNC} We parameterize $V(\vx)$ as ICNN$(6,12,12,1)$, $\alpha(x)$ as $\mathcal{K}(1,10,10,1)$, $\vu(\vx)$ as Control$(6,18,18,6)$. We mask the output of the Control$(6,18,18,6)$ with matrix $[0,0,0,1,1,1]$ s.t. the control only acts on the $\vy$ term. We set $\varepsilon=1\mathrm{e}$-$3$, $c=-0.1$. We train the parameters with $\mathrm{lr}=0.1$ for $300$ steps.

\paragraph{Kernel} We sample $10000$ points from the safe region as the sample from initial distribution $\pi_0$, we set $10000$ zero data as the samples of the target distribution $\pi_1$. We set the  RBF kernel as $k(\vz_1,\vz_2)=\exp(-\|\vz_1-\vz_2\|^2/h),~\vz=(\vx^\top,\vy^\top)^\top$ with bandwidth $h=1e-3$. To move the controlled system from $\pi_0$ to $\pi_1$, the controller is designed as,
\begin{equation*}
\begin{aligned}
    \vu(\vz,t) = \mathbb{E}_{\tilde{\vz}_0\sim\pi_0,\tilde{\vz}_1\sim\pi_1}\left[\dfrac{\tilde{\vz}_1-\vz}{1-t}\dfrac{k(\tilde{\vz(t)},\vz)}{\mathbb{E}_{\tilde{\vz}_0\sim\pi_0,\tilde{\vz}_1\sim\pi_1}[k(\tilde{\vz(t)},\vz)]}\right]-\vf(\vz).
\end{aligned}
\end{equation*}
Here we take the empirical expectation using samples $\tilde{\vz}_0\sim\pi_0$ and $\tilde{\vz}_1\sim\pi_1$, and we denote $\vf(\vz)$ by the original drift term of the 3-link pendulum.

\paragraph{FESS+Kernel} We directly apply our approximate projection operators in Eqs.~\eqref{sa proj}\eqref{st proj} to the Kernel controller, and denote the obtained controller as FESS+Kernel.

We test these controllers and summarize the results in Fig.~\ref{3-link}. In can be seen that although the kernel machine based controller performs poorly in this task, our framework can successfully improve the Kernel controller with safety and stability guarantee.
}

\end{document}